\documentclass[11pt, reqno]{amsart}
\usepackage{amsfonts, amsmath, amssymb, latexsym, eucal, amscd} 

\usepackage{graphicx}
\usepackage{pstricks,pst-plot,pst-node,pst-text,pst-3d}

\begin{document}
 
%%%%%%%%%%%%%%%%%%% 
\numberwithin{equation}{section}

%%%%%%%%%%%%%%%%%%% 
\newtheorem{definition}[equation]{Definition}
\newtheorem{theorem}[equation]{Theorem}
\newtheorem{lemma}[equation]{Lemma}
\newtheorem{corollary}[equation]{Corollary}
\newtheorem{proposition}[equation]{Proposition}
\newtheorem{remark}[equation]{Remark}
\newtheorem{remarks}[equation]{Remarks}
\newtheorem{example}[equation]{Example}
\newtheorem{conjecture}[equation]{Conjecture}
\newtheorem{problem}[equation]{Problem}
\newtheorem{note}[equation]{Note}
 
\def\half{\frac{1}{2}} 
\def\C{\mathbb C}
\def\K{\mathbb K} 
\def\Z{\mathbb Z}
\def\K{\mathbb F}
\def\N{\mathbb N}
\def\I{\mathbb I}
\def\fld{\mathbb K} 

\newcommand{\ad} {\hbox{\rm ad\,}}
\newcommand{\aut}{\hbox{\rm aut}}
\newcommand{\uqhat}{U_q(\widehat{\mathfrak{sl}}_2)}
\newcommand\g {\mathfrak{g}}
\newcommand{\ot}{\otimes}

%%%%%
 \setcounter{secnumdepth}{2}

\title[]{\large The Equitable Basis for $\mathfrak{sl}_2$} 

\author[]{Georgia Benkart$^{\star}$} \address{Department of Mathematics \\ University
of Wisconsin \\  Madison, WI  53706, USA}
\email{benkart@math.wisc.edu}

\author[]{Paul Terwilliger \address{ \\ }}
\email{terwilli@math.wisc.edu}

\thanks{$^{\star}$Support from NSF grant \#{}DMS--0245082 is gratefully acknowledged.  \hfil \break 
{\bf Keywords}:   equitable basis, modular group, Kac-Moody algebra, Cartan matrix     
 \hfil\break
\noindent {\bf 2000 Mathematics Subject Classification}:   17B37}   
 
 \date{January 6, 2010}
  \maketitle
  
 \begin{abstract}  
This article contains an investigation of the equitable basis
for the Lie algebra $\mathfrak{sl}_2$. Denoting this basis by
 $\lbrace x,y,z\rbrace$,   we have
\begin{eqnarray*}
  [x,y] = 2x + 2y,  \qquad [y,z] = 2y + 2z, \qquad [z, x] = 2z + 2x.
\end{eqnarray*}  We determine the
group of automorphisms $G$ generated by
 $\exp(\ad x^*)$, \,$\exp(\ad y^*)$, \,$\exp(\ad z^*)$, where $\lbrace x^*,y^*,z^*\rbrace $ is the
 basis for $\mathfrak{sl}_2$ dual to $\lbrace x,y,z\rbrace $ with respect to the trace form $(u,v) = \hbox{\rm tr}(uv)$ and study the relationship of $G$ to the isometries of
 the lattices $L={\mathbb Z}x \oplus {\mathbb Z}y\oplus {\mathbb Z}z$
 and $L^* ={\mathbb Z}x^* \oplus {\mathbb Z}y^*\oplus {\mathbb Z}z^*$.
 The matrix of the trace form is a Cartan matrix of hyperbolic type, and
 we identify the equitable basis with a set of simple roots 
 of the corresponding Kac-Moody Lie algebra $\mathfrak{g}$, so
 that $L$ is the root lattice and $\half L^*$ is the weight lattice of $\mathfrak g$.  The orbit
 $G(x)$ of $x$ coincides with the set of real roots of $\mathfrak g$. We determine 
 the isotropic roots of $\mathfrak g$ and show that each isotropic
 root has multiplicity 1. We  describe the finite-dimensional $\mathfrak{sl}_2$-modules 
from the point of view of the equitable basis. 
In the final section, we establish a connection between  the Weyl group orbit of the fundamental weights of $\mathfrak{g}$  
and Pythagorean triples.
 \end{abstract} 
\medskip 
 
 \section{
Introduction}   The purpose of this article is to investigate systematically 
a certain basis, called the equitable basis, 
for the Lie algebra $\mathfrak{sl}_2$  of $2 \times 2$ trace zero matrices over
a field $\mathbb F$ of characteristic zero. 
This basis has already appeared in the theory of tridiagonal pairs \cite{H1}, \cite{H2}
and of the three-point loop algebra $\mathfrak{sl}_2 \otimes \mathbb F[t,t^{-1}, (t-1)^{-1}]$
(\cite{HT}, \cite{BT}, see also
\cite{ITW}). 
As we will show, it exhibits many striking features
and has connections with the theory of Kac-Moody Lie algebras.

In Section 2,  we introduce the equitable basis 
 $\lbrace x,y,z\rbrace$  and its dual basis $\lbrace x^*,y^*,z^*\rbrace$
 with respect to the trace form $(u,v) = \hbox{\rm tr}(uv)$ on $\mathfrak{sl}_2$.  
 In Section 3, we study the group $G$ generated by
 $\exp(\ad x^*)$, \,$\exp(\ad y^*)$, \,$\exp(\ad z^*)$ and 
  show in Theorem  \ref{thm1} that $G$ is isomorphic to the modular group
 $\hbox{\rm PSL}_2({\mathbb Z})$.   We then turn our attention to the
 lattices $L={\mathbb Z}x \oplus {\mathbb Z}y\oplus {\mathbb Z}z$
 and $L^* ={\mathbb Z}x^* \oplus {\mathbb Z}y^*\oplus {\mathbb Z}z^*$,
 and in Theorem \ref{thm2},  give a characterization of the orbit $G(x)$ as 
 the elements $u \in L$ with $(u,u) = 2$.  Since there is an automorphism
 of $\mathfrak{sl}_2$ in $G$ which cyclically permutes $x,y,z$,  this is
 the same as the orbit of $y$ and of $z$.   The next two sections are
 devoted to a study of the isometries, automorphisms, and antiautomorphisms
 of $L$ and $L^*$.  
 Theorems  \ref{thm:isostab}-\ref{thm:autL} and \ref{thm:starsame} give the precise relationship
 between (i) the group $G$, \ (ii) the group of automorphisms for $\mathfrak{sl}_2$ that
 preserve $L$,  (iii) the group of automorphisms and antiautomorphisms for $\mathfrak{sl}_2$
 that preserve $L$, and (iv) the group of isometries $\hbox{\rm Isom}_{\mathbb Z}(L)$ for $(\,,\,)$ that preserve
 $L$, as well as analogous results for the lattice $L^* ={\mathbb Z}x^* \oplus {\mathbb Z}y^*\oplus {\mathbb Z}z^*$.  
 
 In Section 7, we make explicit the connections between the equitable basis and the hyperbolic Kac-Moody Lie algebra $\mathfrak g$ corresponding to the Cartan matrix
 $$\left[\begin{array}{ccc} \ \ 2 & -2 & -2 \\ -2 & \ \ 2 & -2 \\ -2 & -2 &\ \ 2 \end{array} \right].$$
 The simple roots of $\mathfrak g$ can be identified  with the elements of the
 equitable basis, and  the real roots with the orbit $G(x)$ (see Theorem \ref{thm3}).   
 In Proposition \ref{prop:iso2} we show that the set of  isotropic roots can be
 identified with $\bigcup_{n \in \mathbb Z, n \neq 0}    2nG (z^*)$, which is
 precisely the set of nonzero nilpotent matrices in $\mathfrak{sl}_2(\mathbb Z)$,
 and prove that each  isotropic root has multiplicity 1 (see Corollary 
 \ref{cor:isomult}).   We determine the relationship between the Weyl
 group $W$  of $\mathfrak g$ and the group $\hbox{\rm Isom}_{\mathbb Z}(L)$
 in Proposition \ref{prop:isoWeyl}. 
 
 Section 8 studies the finite-dimensional representations of $\mathfrak{sl}_2$
 from the equitable point of view.  Then starting with the equitable
 picture of the adjoint representation for
 $\mathfrak{sl}_2$, in the final section we apply reflections in
 $W$ to obtain the Poincar\'e disk.  
 The equitable basis enables us to connect the Weyl group orbit of the fundamental weights of $\mathfrak{g}$
 with Pythagorean triples.

\medskip

 \section{The equitable basis} 

Throughout,   $\{e,f,h\}$ will denote the basis for $\mathfrak{sl}_2$
given by 

$$  
e = \left[\begin{array}{cc}0 & 1 \\ 0 & 0 \end{array}\right],  \quad  
f = \left[\begin{array}{cc}0 & 0 \\ 1 & 0 \end{array}\right],  \quad
h = \left[\begin{array}{cc}1 &\ \ 0 \\ 0 & -1 \end{array}\right],  $$

\noindent and having products  \ $[e,f] = h, \ [h,e] = 2e, \  [h,f] = -2f$.   \medskip

The {\it equitable basis} $\{x,y,z\}$ for $\mathfrak{sl}_2$  consists of the matrices

\begin{eqnarray}
x &=& \left[\begin{array}{cc}1 & \ \  0 \\ 0 & -1 \end{array}\right] = h,  \nonumber \\  
y &=& \left[\begin{array}{cc} -1 & 2 \\  \ \ 0  & 1 \end{array}\right] = 2e - h, \label{eqbas}  \\
z &=& \left[\begin{array}{cc} -1 &  0 \\  -2  & 1 \end{array}\right]  = -2f - h, \nonumber
\end{eqnarray}

\noindent whose products satisfy

\begin{equation}\label{eq:mult}  [x,y] = 2x + 2y,  \qquad [y,z] = 2y + 2z, \qquad [z, x] = 2z + 2x. \end{equation}
{F}rom this it follows that there is a Lie algebra automorphism $\varrho$ of
$\mathfrak{sl}_2$ of order $3$  such  that

\begin{equation}\label{eq:rho}  \varrho(x) = y, \qquad \varrho(y) = z, \qquad \varrho(z) = x. \end{equation}
Note that
$$ y = \exp(\ad e)(-h), \qquad \qquad z = \exp(\ad f)(-h),$$
where  $\ad u (v) = [u,v]$ and
$\exp (w) =\sum_{n=0}^\infty w^n/n!$. 
We will relate the automorphisms 
$\exp(\ad e)$ and $\exp(\ad f)$ to $\varrho$ in Section 3.

\medskip In our work we will use  the trace form $(u,v) := \hbox{\rm tr}(uv)$ for $u,v \in \mathfrak{sl}_2$.
 We could use instead the Killing form 
$\kappa(u,v) := \hbox{\rm tr}(\ad u\, \ad v)  = 4(u,v)$, but
the trace form has some aesthetic advantages.        
Relative to the equitable basis,  the matrix of the trace form is
given by

\begin{equation}\label{eq:CM}\mathcal A  = \left[\begin{array}{ccc} \ \ 2 & -2 & -2 \\ -2 & \ \ 2 & -2 \\ -2 & -2 &\ \ 2 \end{array} \right].
\end{equation}
This is a (generalized) Cartan matrix as defined in (\cite[\S1.1]{K}, \cite[\S 3.4]{MP}); the
corresponding Kac-Moody Lie algebra will be related to the
equitable basis in Section 6. 

\medskip  Let $\{x^*, y^*, z^*\}$ denote the basis for $\mathfrak{sl}_2$ that is dual 
to the equitable basis in the sense that
$(u,v^*) = 2 \delta_{u,v}$ for all $u,v \in \{x,y,z\}$  (the factor of $2$ is inessential but convenient).
Then
\begin{equation} x+y=-2z^*, \qquad 
y+z=-2x^*,\qquad z+x=-2y^*
\label{eq:1p5}
\end{equation}
and

\begin{equation}  \varrho(x^*) = y^*, \qquad \varrho(y^*) = z^*, \qquad \varrho(z^*) = x^*. \end{equation}
Relative to the basis $\lbrace x^*,y^*,z^*\rbrace $ the matrix of the trace form is 
\begin{equation}\label{eq:CMinv} 4\mathcal A^{-1}  = \left[\begin{array}{ccc} \ \ 0 & -1 & -1 \\ -1 & \ \ 0 & -1 \\ -1 & -1 &\ \ 0 \end{array} \right].
\end{equation}
The equitable basis and its dual are related by the following
multiplication tables: 

\begin{equation}\label{eq:tab1-2}
\begin{tabular}[t]{|c||c|c|c|}
\hline
        $[\, ,\,]$  
       & $x^*$ & $y^*$ & $z^*$
	\\
	\hline  \hline
$x^*$ & 
$0$ & $z$ & $-y$ 
\\  \hline
$y^*$ & 
$-z$ & $0$ & $x$   
\\  \hline
$z^*$ & 
$y$ & $-x$ & $0$  \\
\hline
\end{tabular}
\qquad \qquad
\begin{tabular}[t]{|c||c|c|c|}
\hline
        $[\, ,\,]$  
       & $x$ & $y$ & $z$
	\\
	\hline  \hline
$x$ & 
$0$ & $-4z^*$ & $4y^*$ 
\\  \hline
$y$ & 
$4z^*$ & $0$ & $-4x^*$   
\\  \hline
$z$ & 
$-4y^*$ & $4x^*$ & $0$  \\
\hline
\end{tabular}
\end{equation}

We also have
\begin{equation}
\label{tablexxs}
\begin{tabular}[t]{|c||c|c|c|}
\hline 
        $[\, ,\,]$  
       & $x$ & $y$ & $z$
	\\
	\hline  \hline
$x^*$ & 
$y-z$ & $y+z$ & $-y-z$ 
\\  \hline
$y^*$ & 
$-z-x$ & $z-x$ & $z+x$ 
\\  \hline
$z^*$ & 
$x+y$ & $-x-y$ & $x-y$ \\  \hline
\end{tabular}  \qquad  \quad
\end{equation}
and 
\begin{eqnarray}
&& x-y=2(x^*-y^*), \qquad 
y-z=2(y^*-z^*), \qquad
z-x=2(z^*-x^*), \nonumber
\\
&& x^*+y^*=z^*-z,\qquad
y^*+z^*=x^*-x,\qquad
z^*+x^*=y^*-y, \label{eq:3sum}
\\
&& \qquad \qquad \qquad \qquad x+y+z=-x^*-y^*-z^*.\nonumber
\end{eqnarray}
\medskip
By (\ref{eqbas}),  each of the matrices $x,y,z$ is semisimple (diagonalizable) with eigenvalues
$1$ and $-1$. Since

\begin{equation} x^* = h-e+f, \qquad  y^* = f, \qquad z^* = -e, 
\label{eq:shownil} \end{equation}

\noindent each of the dual basis elements $u\in \{x^*, y^*, z^*\}$  is nilpotent with  $u^2 = 0$
$(\ad u)^3 = 0$.   
\medskip

\section{Connections with the modular group}

Let $G$ denote the subgroup of the automorphism group  $\hbox{\rm Aut}_{\mathbb F}(\mathfrak {sl}_2)$
generated by $\exp(\ad x^*)$, \,$\exp(\ad y^*)$, \,$\exp(\ad z^*)$.  
In this section we will prove that $G$ is isomorphic to the modular group $\hbox{\rm PSL}_2(\mathbb Z)$.   Recall that $\hbox{\rm PSL}_2(\mathbb Z)$ is obtained from the
group $\hbox{\rm SL}_2(\mathbb Z)$ of $2 \times 2$ integral matrices of determinant 1
by factoring out the subgroup consisting of the matrices $\pm I$.     
It  is a free product of a cyclic group of order 2 and a cyclic group of order 3 (see for example, \cite{A1}).   
To establish the isomorphism with $G$,  we first locate generators
for $G$ of order 2 and 3.    \bigskip

\begin{definition} \label{def1}  Let $\sigma_x$, $\sigma_y$, $\sigma_z$ be the automorphisms
of $\mathfrak{sl}_2$ defined by
$$ \sigma_x = \exp(\ad x^*), \qquad \sigma_y = \exp(\ad y^*),  \qquad \sigma_z = \exp(\ad z^*).$$
\end{definition}
\bigskip

Using the table in \eqref{tablexxs} we obtain \bigskip

\begin{lemma}\label{lem0}  The matrices of $\sigma_x,\sigma_y,\sigma_z$ relative
to the equitable basis are given by

\begin{equation*}
\sigma_x \rightarrow \left[\begin{array}{ccc} 1 & 0 & \ \ 0 \\ 2 & 2 & -1 \\ 0 & 1 & \ \ 0 \end{array}\right],
\qquad \sigma_y \rightarrow \left[\begin{array}{ccc} \ \ 0 & 0 & 1 \\ \ \ 0 & 1 & 0 \\ -1 & 2 & 2 \end{array}\right], \qquad
\sigma_z \rightarrow \left[\begin{array}{ccc} 2 & -1 & 2 \\ 1 & \ \ 0 & 0 \\ 0 & \ \ 0 & 1 \end{array}\right].
\end{equation*}
\end{lemma}
\medskip

\begin{lemma}\label{lem1}
\begin{itemize}  \item[(i)]  $\varrho \sigma_x \varrho^{-1} = \sigma_y$,  \qquad $\varrho \sigma_y \varrho^{-1} = \sigma_z$,
 \qquad $\varrho \sigma_z  \varrho^{-1} = \sigma_x$;

\item [(ii)] $\varrho$ is equal to each of the products $\sigma_x \sigma_y$, 
$\sigma_y \sigma_z$, $\sigma_z \sigma_x$.   In particular $\varrho \in G$.
\end{itemize} 
\end{lemma}

\noindent {\it Proof:}  Part (i)  follows from the well-known identity
$$\varphi \exp(\ad u) \varphi^{-1} = \exp(\ad \varphi(u))$$
which holds for all $\varphi \in \hbox{\rm Aut}_{\mathbb F}(\mathfrak{sl}_2)$ and  all
nilpotent $u \in \mathfrak {sl}_2$.  
To see that $\varrho= \sigma_x\sigma_y$,   use 
\eqref{eq:rho}
and Lemma \ref{lem0} to 
verify that $\varrho$ and $\sigma_x \sigma_y$
agree on  the elements of the equitable basis. To obtain the other two
expressions for $\varrho$,   apply  part (i) above. \qed
 
\bigskip
\begin{lemma}\label{lem2}  For $\sigma_x$, $\sigma_y$, $\sigma_z$ as in
Definition \ref{def1},  we
have the  following: 
\begin{itemize}
\item[(a)]   $\sigma_x = \sigma_y \sigma_z \sigma_y^{-1}$;
\item[(b)]   $\sigma_x = \sigma_z^{-1} \sigma_y \sigma_z$;
\item[(c)]   $(\sigma_y \sigma_z)^3 = 1$;
\item[(d)]   $\sigma_y \sigma_z \sigma_y = \sigma_z \sigma_y \sigma_z$;
\item[(e)]   $(\sigma_y \sigma_z \sigma_y)^2 = 1$;
\item[(f)]    $G$ is generated by $\sigma_y$ and $\sigma_z$.  
\end{itemize}
\end{lemma}

\noindent {\it Proof:}   These properties can be deduced from Lemma \ref{lem1}.
\qed 

\bigskip     

\begin{definition}\label{def2}  Let  $\tau_x$, $\tau_y$, $\tau_z$ denote the elements of $G$
defined by 
\begin{eqnarray*} \tau_x &=& \sigma_y \sigma_z \sigma_y = \exp(\ad y^*) \exp(\ad z^*) \exp(\ad y^*) \\
&=& \sigma_z \sigma_y \sigma_z =  \exp(\ad z^*) \exp(\ad y^*)  \exp (\ad z^*), \\
\tau_y &=&  \sigma_z \sigma_x \sigma_z = \exp(\ad z^*) \exp(\ad x^*) \exp(\ad z^*) \\
&=& \sigma_x \sigma_z \sigma_x = \exp(\ad x^*) \exp(\ad z^*)  \exp (\ad x^*), \\
\tau_z &=&  \sigma_x \sigma_y \sigma_x = \exp(\ad x^*) \exp(\ad y^*) \exp(\ad x^*) \\
&=&  \sigma_y \sigma_x \sigma_y   = \exp(\ad y^*) \exp(\ad x^*)  \exp (\ad y^*). 
\end{eqnarray*}

\end{definition}

\medskip
\begin{lemma}\label{lem3}  
For $\tau_x$, $\tau_y$, $\tau_z$ as in Definition \ref{def2},   the following relations hold:
\begin{itemize}
\item[(a)]  $\varrho \tau_x \varrho^{-1} = \tau_y,  \qquad
\varrho \tau_y \varrho^{-1} = \tau_z, \qquad  
\varrho \tau_z \varrho^{-1} = \tau_x$;
\item[(b)]  $\tau_x^2 = \tau_y^2 = \tau_z^2 = 1$;
\item[(c)]  $\sigma_z = \tau_x \varrho^{-1}$ and $\sigma_y = \varrho^{-1} \tau_x$;
\item[(d)]  $\tau_x(x) = -x$,  \ $\tau_x(y) = 2x + z$,  \  $\tau_x(z) = 2x+y$.  
\end{itemize}
\end{lemma} 

\noindent {\it Proof:}   Part (a) follows from Lemma \ref{lem1}\,(i),
while (b) comes from (a) and Lemma \ref{lem2}\,(e).
Concerning (c), the first (resp. second) equation follows from $\varrho=\sigma_y\sigma_z$
and $\tau_x=\sigma_z\sigma_y\sigma_z$
(resp. $\tau_x=\sigma_y\sigma_z\sigma_y$).
To get (d),  use $\tau_x=\sigma_z \varrho$ together with \eqref{eq:rho} and Lemma \ref{lem0}.
\qed 
\medskip

Combining Lemma \ref{lem3} with Lemma \ref{lem2}\,(f),  we have 

\bigskip
\begin{corollary}\label{cor1}  Each of the following is a generating set
for the group $G$. 
\smallskip

$\hbox{\rm (i)} \ \varrho,\, \tau_x$;   \qquad  \hbox{\rm (ii)}  $\varrho, \, \tau_y$;  \qquad \hbox{\rm (iii) }$\varrho, \tau_z$.

\end{corollary}

\medskip
For  $\theta \in \hbox{\rm SL}_2(\mathbb Z)$, conjugation by $\theta$
determines an automorphism $\widehat \theta$ of $\mathfrak{sl}_2$:

$$\widehat \theta:   u \mapsto \theta u \theta^{-1}.$$

\noindent The map   

\begin{equation}\label{eq:thetahat}  \begin{array}{ccc}    \hbox{\rm SL}_2(\mathbb Z) &\rightarrow& \hbox{\rm Aut}_{\mathbb F}(\mathfrak{sl}_2) \\
\theta &\mapsto&  \widehat \theta \end{array} \end{equation}
is a group homomorphism with kernel $\{\pm I\}$. This map induces an embedding  

\begin{equation}   \imath:  \hbox{\rm PSL}_2(\mathbb Z)  \rightarrow \hbox{\rm Aut}_{\mathbb F}(\mathfrak{sl}_2). 
\end{equation}

\bigskip
\begin{theorem} \label{thm1}  The image of the embedding
$ \imath:  \hbox{\rm PSL}_2(\mathbb Z)  \rightarrow \hbox{\rm Aut}_{\mathbb F}(\mathfrak{sl}_2) $
coincides with $G$. 
Therefore 
$G$ is isomorphic to  $\hbox{\rm PSL}_2(\mathbb Z)$.   
\end{theorem}

\noindent {\it Proof:} The matrices 
\begin{equation} A=\left [\begin{array}{cc}  0& -1 \\ 1 & \ 1 \end{array} \right],  \qquad 
B=\left [\begin{array}{cc}  0 & -1 \\ 1 & \ 0 \end{array} \right],
\qquad C=\left [\begin{array}{cc}  1& -1 \\ 1 & \ 0 \end{array} \right]
\end{equation} 
are in $\hbox{\rm SL}_2(\mathbb Z)$
and satisfy $C=BAB^{-1}$. 
Let $a, b, c$ denote the images of $A,B,C$ respectively under the canonical
homomorphism $\hbox{\rm SL}_2(\mathbb Z) \to
\hbox{\rm PSL}_2(\mathbb Z)$, 
and note that $c=bab^{-1}$. 
By \cite{A2},  the elements $a,b$ generate
$\hbox{\rm PSL}_2(\mathbb Z)$,  so $b,c$ generate $\hbox{\rm PSL}_2(\mathbb Z)$.
One checks that ${\widehat B}=\tau_x$ and
 ${\widehat C}=\varrho$ so
$\imath (b)=\tau_x$ and
$\imath (c)=\varrho$.
The result then follows in view of Corollary \ref{cor1}\,(i). \qed

\medskip
\section{The $G$-orbit of $x$}

In this section we describe the orbit $G(x)$  of $x$ under
the group $G $ generated by $\exp(\ad x^*), \, \exp(\ad y^*), \, \exp(\ad z^*)$.  
Since $\varrho$ belongs to $G$ and cyclically permutes
the elements of the equitable basis, $G(x)$ coincides with the
$G$-orbit of $y$ and the $G$-orbit of $z$.   Later in the paper we relate $G(x)$
to the set of real roots for the Kac-Moody Lie algebra associated with the Cartan matrix 
$\mathcal A$ from
\eqref{eq:CM}.    We begin by determining the stabilizer of $x$ in $G$.  

\bigskip
\begin{lemma}\label{lem4}  Suppose $g \in G$ and $g(x) = x$.   Then $g = 1$. 
\end{lemma}

\noindent {\it Proof:}  
By Theorem \ref{thm1} and the paragraph preceding it, there exists
$\theta \in \hbox{\rm SL}_2(\mathbb Z)$ such that ${\widehat \theta}=g$. Therefore
 $\theta x \theta^{-1} = g(x)=x$  gives  $\theta x = x \theta$, and  this along with the fact that $x=\mbox{\rm diag}(1,-1)$ implies 
 $\theta$ is diagonal. The diagonal entries of $\theta$ are integers whose product is
 1,  so they are both 1 or both $-1$; thus $\theta=\pm I$ so $g=1$. \qed

\medskip
\begin{corollary}\label{cor2}  The map 
\begin{equation*} \begin{array}{ccc}  G &\rightarrow& G(x) \\
g& \mapsto& g(x)  \end{array}   \end{equation*}
is a bijection.  \end{corollary}

\bigskip
We turn our attention now  to the lattice 
\begin{eqnarray}\label{eq:lattice}
L: = \mathbb Z x \oplus  \mathbb Z y  \oplus \mathbb Z z.
\end{eqnarray}  Our  goal is to prove that

$$ G(x) = \{u \in L \mid (u,u) = 2\}, $$

\noindent but this necessitates a few comments about $L$.   Note that $L$ is closed under the Lie
bracket and invariant under the group $G$.    Further observe that

\begin{equation}\label{eq:mod2}  L = \left \{ \left [\begin{array}{cc}  p & \ \ q \\ r & -p \end{array} \right]  
\, \Bigg | \,p,q,r \in \Z,\quad   q,r \; {\rm even} \right \}.  \end{equation}

\noindent This realization
of $L$ shows that it is
the Lie algebra analogue of the congruence subgroup of $\hbox{\rm PSL}_2(\mathbb Z)$.
\medskip
 
Recall that an {\it isometry} of $\mathfrak{sl}_2$ is an $\K$-linear bijection
$\varphi:\mathfrak{sl}_2 \to \mathfrak{sl}_2$ such that $(\varphi(u),\varphi(v))=(u,v)$
for all $u,v\in \mathfrak{sl}_2$.
\medskip

\begin{lemma}\label{eq:ginv}
Each automorphism of $\mathfrak{sl}_2$ is an isometry of $\mathfrak{sl}_2$.
\end{lemma}
\noindent {\it Proof:}  Each automorphism of a finite-dimensional Lie algebra is
an isometry relative to the Killing form.    Since the trace map is a multiple
of the Killing form, the result is apparent.  \qed 
 
\medskip

The next  lemma provides some useful formulae for the square norm $(u,u)$ of an
element $u \in \mathfrak{sl}_2$.  \bigskip

\begin{lemma} \label{lem:4a}
For $u = \alpha x +\beta y+\gamma z \in \mathfrak{sl}_2$,  the expression $(u,u)/2$
is equal to each of the following:
\begin{eqnarray*}
&&\alpha^2 + \beta^2 +\gamma^2 -2 (\alpha\beta + \beta\gamma +\gamma\alpha),
\\
&&2\big(\alpha^2 +\beta^2 +\gamma^2\big) - (\alpha+\beta+\gamma)^2,
\\
&& \big(\alpha +\beta -\gamma\big)^2 - 4\alpha \beta,
\\
&& \big(\beta +\gamma -\alpha \big)^2 - 4\beta \gamma, 
\\
&& \big(\gamma+\alpha-\beta\big)^2 - 4\gamma \alpha.
\end{eqnarray*}
\end{lemma}
\noindent {\it Proof:} The above five scalars are mutually equal; this can
be checked by algebraic manipulation. 
Observe that $(u,u)=(\alpha,\beta,\gamma){\mathcal A}(\alpha,\beta,\gamma)^t$
where $\mathcal A$ is from \eqref{eq:CM}. Evaluating this triple product by matrix
multiplication,  we find that $(u,u)/2$ is equal to the first expression above, so
the result follows. \qed

\begin{definition} 
\label{def:r}
\rm Define  $R = \{u \in L \mid (u,u) = 2\}$.   We note
that $R$ is $G$-invariant by Lemma \ref{eq:ginv}.
\end{definition} 
\medskip

\begin{lemma}\label{lem5}
For $u=\alpha x +\beta y +\gamma z$ in the set $R$ in Definition \ref{def:r},  
 the coefficients
$\alpha,\beta,\gamma$ are either all nonnegative or all are nonpositive.   \end{lemma}

\noindent {\it Proof:}  
We assume the result is false and reach a contradiction. There exists a pair of coefficients
having  opposite signs; without loss in generality we may assume they are $\alpha$ and $\beta$.
Thus $\alpha \beta \leq -1$. Using $(u,u)=2$ and the third expression in Lemma
\ref{lem:4a},  we find that 
$$-4 \geq 4 \alpha \beta  = (\alpha + \beta -\gamma)^2 - 1  \geq  -1.$$
This is a contradiction,  so the result must be true.   \qed 
 
\bigskip

\begin{definition} \rm For the set $R$ in Definition \ref{def:r} 
and for $u = \alpha x + \beta y + \gamma z \in R$,  we define
the {\it height} of $u$ to be the sum 
$\hbox{\rm ht}(u) = \alpha+\beta+\gamma$.   Let
$R^+ = \{ u \in R \mid  \hbox{\rm ht}(u) > 0\}$ and $R^- = \{ u \in R \mid  \hbox{\rm ht}(u) < 0\}$.  
\end{definition}

\medskip

By Definition 
\ref{def:r} and
 Lemma \ref{lem5} we have $R = R^+ \cup R^-,$  $R^- = -R^+$,  
and $\varrho (R^\pm)  = R^\pm$. Next we describe how the automorphisms
$\tau_x$, $\tau_y$, $\tau_z$ act on the set $R^+$.  
\bigskip

\begin{lemma}\label{lem6} For $u \in \lbrace x,y,z\rbrace $, 
the map $\tau_u$ sends $u$ to $-u$ and permutes the elements of 
$R^+ \setminus \{u\}$.   \end{lemma}
 
\noindent {\it Proof:} There is no loss in generality in assuming  $u=x$.
Recall that $\tau_x(x) = -x$ by
Lemma \ref{lem3}\,(d). Now suppose we are given $v  = \alpha x + \beta y + \gamma z
\in R^+$ such that $\tau_x(v) \in R^-$.  It suffices to argue
that $v = x$.  Using Lemma
\ref{lem3}(d), we have 

\begin{eqnarray*}  \tau_x(v) &=& \alpha(-x) + \beta(2x+z) + \gamma(2x+y)  \\
&=& (2\beta+2\gamma-\alpha)x+ \gamma y + \beta z. \end{eqnarray*}
Observe $\beta\geq 0$ since $v \in R^+$ and $\beta \leq 0$ since
$\tau_x(v) \in R^-$ so $\beta=0$. Similarly $\gamma=0$.
Now $\alpha=1$ since $(v,v) = 2$. Therefore $v=x$ and the result follows.
\qed

\medskip
\begin{theorem}  \label{thm2}  $G(x) = R = \{u \in L \mid (u,u) = 2\}$.   \end{theorem}

\noindent {\it Proof:} The set $R$ contains $x$ and is $G$-invariant so $G(x)\subseteq R$.
To show equality holds, we assume there exists $u \in R \setminus G(x)$ and arrive at a contradiction.
Without loss we may further assume that $u \in R^+$ and that
$u$ has minimal height with this property. Note that $u$ is not one of $x,y,z$
as they are in $G(x)$. 
Write
$u=\alpha x + \beta y + \gamma z$.  By Lemma \ref{lem6}, each of 
$\tau_x(u)$, $\tau_y(u)$, $\tau_z(u)$ belongs to $R^+ \setminus G(x)$. By our minimality assumption
all these elements have height at least $\hbox{\rm ht}(u)$. Evaluating these
inequalities we determine that 
\begin{eqnarray*} 
\alpha \leq  \beta +\gamma, \qquad 
\beta\leq \gamma +\alpha, \qquad 
\gamma\leq \alpha + \beta.
\end{eqnarray*}
Since the situation is cyclically symmetric,  
we may assume that $\alpha \geq \beta$.   Using
$(u,u)=2$ and the third expression in 
Lemma \ref{lem:4a} 
we see that  
\begin{equation*}  (\alpha +\beta-\gamma)^2 = 1+4\alpha \beta  > 4 \beta^2, \end{equation*}
so that $\alpha + \beta -\gamma > 2 \beta$ and then $\alpha > \beta + \gamma$.  
This is a contradiction,  so it must be that  $G(x)=R$.  \qed
\medskip
 
\begin{note}\label{cor3} \rm
Combining Theorem \ref{thm1}, Corollary \ref{cor2}, and Theorem \ref{thm2} and using
Theorem \ref{lem:4a},
we get a bijection between $\hbox{\rm PSL}_2(\mathbb Z)$ and the set of
integral solutions $(\alpha,\beta,\gamma)$ to the quadratic equation

$$2\big(\alpha^2 +\beta^2 +\gamma^2\big) - (\alpha+\beta+\gamma)^2 = 1.$$

\end{note}
\medskip

We close this section with a result about the coefficients of elements of $G(x)$.

\bigskip
\begin{proposition}\label{prop4}   For $\alpha x + \beta y + \gamma z \in R$, exactly one
of the coefficients $\alpha,\beta,\gamma$ is odd.  
\end{proposition}

\noindent {\it Proof:}   Let $S$ denote the set of  elements in $G(x)$ with exactly one odd coefficient.
Note that $S$ contains $x$.    For $u =\alpha x + \beta y + \gamma z \in S$ we have
$\tau_x(u) = (2\beta+2\gamma-\alpha)x + \gamma y + \beta z$, and modulo 2,  this
element has the same coefficients as $u$.  Therefore $\tau_x(u) \in S$.     Since $\varrho(u) \in S$
also, we see that $S$ is $G$-invariant.   Consequently  $G(x) \subseteq S$ so $G(x)=S$.   \qed

\medskip
\section{Automorphisms, antiautomorphisms, and isometries}

In this section we continue our study of the lattice
$L$ from \eqref{eq:lattice}.  We  use the equitable basis 
to determine  the precise relationship
between the following four groups: 
(i) the group $G$ from Section 3,  (ii) the group of automorphisms of $\mathfrak{sl}_2$
that preserve $L$,  (iii) the group of automorphisms and antiautomorphisms
of $\mathfrak{sl}_2$ that preserve $L$, and  (iv) the group of isometries of 
$\mathfrak{sl}_2$ that preserve $L$.
 
\medskip
By an {\it antiautomorphism} of $\mathfrak{sl}_2$
we mean an $\K$-linear bijection $\phi: \mathfrak{sl}_2\to \mathfrak{sl}_2$ such
that $\phi(\lbrack u,v\rbrack)= \lbrack \phi(v),\phi(u)\rbrack$ for $u,v \in \mathfrak{sl}_2$.
Here are some examples. The map $-1:u\to -u$ is an antiautomorphism of $\mathfrak{sl}_2$ (in
fact of any Lie algebra).
For distinct $u,v \in \lbrace x,y,z\rbrace$,  the $\K$-linear map   $(u\,v):\mathfrak{sl}_2 \to
\mathfrak{sl}_2$  that interchanges $u,v$ and fixes the remaining element in $\lbrace x,y,z\rbrace $
  is an antiautomorphism of $\mathfrak{sl}_2$.
Let $\hbox{\rm AAut}_{\mathbb F}(\mathfrak {sl}_2)$
denote the group consisting of the automorphisms and antiautomorphisms of $\mathfrak {sl}_2$.
Then $\hbox{\rm Aut}_{\mathbb F}(\mathfrak {sl}_2)$ is a normal subgroup of
$\hbox{\rm AAut}_{\mathbb F}(\mathfrak {sl}_2)$ of index 2,  and
\begin{eqnarray}   \hbox{\rm AAut}_{\mathbb F}(\mathfrak{sl}_2)
 = \{\pm 1\} \ltimes \hbox{\rm Aut}_{\mathbb F}(\mathfrak{sl}_2).
\label{eq:aainfo}
 \end{eqnarray}
Define
\begin{eqnarray*}\label{eq:autz} \hbox{\rm AAut}_{\mathbb Z}(L) &=& \{ \varphi \in   \hbox{\rm AAut}_{\mathbb F}(\mathfrak {sl}_2) \mid   \varphi(L) = L\}, \\
\hbox{\rm Aut}_{\mathbb Z}(L) &=& \{ \varphi \in   \hbox{\rm Aut}_{\mathbb F}(\mathfrak {sl}_2) \mid   \varphi(L) = L\}
 \end{eqnarray*} 
and note that
\begin{eqnarray} \label{eq:normal}  \hbox{\rm AAut}_{\mathbb Z}(L)
 = \{\pm 1\} \ltimes \hbox{\rm Aut}_{\mathbb Z}(L).
 \end{eqnarray}   We remark that  $\hbox{\rm AAut}_{\mathbb Z}(L)$  
is the group  of automorphisms and antiautomorphisms of $L$,  viewed as a Lie algebra over $\mathbb Z$,  since every such map over $\mathbb Z$   
 can be extended linearly to an automorphism or antiautomorphism  in $\hbox{\rm AAut}_{\mathbb F}(\mathfrak {sl}_2)$.
 By the construction $G \subseteq \hbox{\rm Aut}_{\mathbb Z}(L)$. 
 Let $\mbox{\rm Isom}_{\K}(\mathfrak{sl}_2)$ denote the group of all isometries
of $\mathfrak{sl}_2$ and define
\begin{eqnarray*}
  \mbox{\rm Isom}_{\Z}(L) = \lbrace \varphi \in 
   \mbox{\rm Isom}_{\K}(\mathfrak{sl}_2) \;|\;\varphi(L)=L \rbrace.
   \end{eqnarray*}
 Since $-1$ is an isometry of $\mathfrak{sl}_2$,  by Lemma  \ref{eq:ginv} 
  we have 
 $\hbox{\rm AAut}_{\mathbb Z}(L)\subseteq \hbox{Isom}_{\mathbb Z}(L)$.
So far we know that 
\begin{eqnarray*}
G \subseteq 
\hbox{\rm Aut}_{\mathbb Z}(L)
 \subseteq 
\hbox{\rm AAut}_{\mathbb Z}(L) 
 \subseteq 
\hbox{\rm Isom}_{\mathbb Z}(L).
\label{eq:chain}
\end{eqnarray*}
Before describing  this chain in more detail 
we compute the stabilizer of $x$ in $\hbox{\rm Isom}_{\mathbb Z}(L)$.
 In what follows $\langle S \rangle $ means the group generated by the set $S$. \bigskip
 
 \begin{lemma} 
\label{lem:isostab}
The stabilizer of $x$ in $\hbox{\rm Isom}_{\mathbb Z}(L)$
 is $\langle (y\,z), -\tau_x\rangle $ where $\tau_x$ is from Definition
 \ref{def2}. Hence, this stabilizer is isomorphic to ${\mathbb Z}_2 \times {\mathbb Z}_2$.
 \end{lemma}
\noindent {\it Proof:}  Concerning the first assertion,
one inclusion is clear since each of the maps $(y\,z)$, $-\tau_x$ fixes $x$ and
is contained in $\hbox{\rm Isom}_{\mathbb Z}(L)$.
To obtain the other inclusion, we pick $\varphi \in 
\hbox{\rm Isom}_{\mathbb Z}(L)$ such that $\varphi(x)=x$ and show $\varphi \in 
\langle (y\,z), -\tau_x\rangle $. 
The subspace $\hbox{\rm Span}_{\mathbb F}\lbrace y^*,z^*\rbrace $ 
is the orthogonal complement of $x$
relative to the trace form, so this subspace is $\varphi$-invariant.
Write $\varphi(y^*)=ay^*+bz^*$ and $\varphi(z^*)=cy^*+dz^*$. Using
$\varphi(x)=x$ and $x+z=-2y^*$, $x+y=-2z^*$, we have that   
$x+\varphi(z)=a(x+z)+b(x+y)$; this shows $a,b \in {\mathbb Z}$
since $\varphi(z) \in L$.     Similarly $c,d \in {\mathbb Z}$.    By 
\eqref{eq:CMinv} the matrix representing  the trace form relative
to $\lbrace y^*, z^*\rbrace $ is
$\left [\begin{array}{cc}  0& -1 \\ -1 & \ 0 \end{array} \right]$.
Since $\varphi$ is an isometry,
\begin{eqnarray*}
\left [\begin{array}{cc}  a& b \\ c & d \end{array} \right] 
\left [\begin{array}{cc}  0& -1 \\ -1 & 0 \end{array} \right]
\left [\begin{array}{cc}  a& c \\ b &  d \end{array} \right] =
\left [\begin{array}{cc}  0& -1 \\ -1 &  0 \end{array} \right]
\end{eqnarray*}
and this yields $ab=0$, $cd=0$, $ad+bc=1$ after a brief calculation. By these equations
 $\left [\begin{array}{cc}  a& c \\ b & d \end{array} \right]$ is one of
 \begin{eqnarray*}
 \left [\begin{array}{cc}  1& 0 \\ 0 &  1 \end{array} \right], \qquad 
\left [\begin{array}{cc}  0& 1 \\ 1 &  0 \end{array} \right], \qquad 
\left [\begin{array}{cc}  0& -1 \\ -1 &  0 \end{array} \right],\qquad
\left [\begin{array}{cc}  -1& 0 \\ 0&  -1 \end{array} \right]
 \end{eqnarray*}
 and these solutions correspond to $\varphi=1$, $\varphi=(y\,z)$, $\varphi=-\tau_x$,
$\varphi=-(y\,z)\tau_x$ respectively. In any event,  $\varphi \in 
\langle (y\,z), -\tau_x\rangle $ and the first assertion follows. 
The second assertion is a direct consequence of the first.
\qed  
\bigskip

\begin{theorem}
\label{thm:isostab}
 $\hbox{\rm Isom}_{\mathbb Z}(L)$ is equal to each of the groups
\begin{eqnarray*}
  \langle  (x\,y),-1 \rangle \ltimes G,
  \qquad 
 \langle (y\,z), -1 \rangle \ltimes G,
  \qquad  \langle (z\,x), -1 \rangle \ltimes G.
\end{eqnarray*} In particular $\hbox{\rm Isom}_{\mathbb Z}(L)$ is isomorphic
to $({\mathbb Z}_2 \times {\mathbb Z}_2) \ltimes G$.
\end{theorem}
\noindent {\it Proof:} 
 We first show that $\hbox{\rm Isom}_{\mathbb Z}(L)=\langle (y\,z), -1 \rangle \ltimes G$.
  Certainly $-1$ normalizes $G$ since $-1$ commutes with everything in $G$.
  The element $(y\,z)$ also normalizes $G$ since
  $(y\,z)\tau_x =\tau_x (y\,z)$,  $(y\,z) \varrho = \varrho^2 (y\,z)$, and  
  $\tau_x, \varrho$ together generate $G$.
  The group $G$ has trivial intersection with
  $\langle (y\,z), -1 \rangle$ because $G$ has trivial intersection
with  $\langle (y\,z), -\tau_x \rangle$ by Lemma \ref{lem4} and Lemma
\ref{lem:isostab}, and because $\tau_x \in G$.
To see that $\hbox{\rm Isom}_{\mathbb Z}(L)$
is generated by $G$, $(y\,z)$, $-1$,   
choose $\varphi \in  \hbox{\rm Isom}_{\mathbb Z}(L)$. Since $\varphi$ is
an isometry of $\mathfrak{sl}_2$,  the set $R$ from Definition \ref{def:r}
is $\varphi$-invariant. Recall that $x \in R$,  so $\varphi(x) \in R$.
But $R=G(x)$ so there exists $g \in G$ such that
$\varphi(x)=g(x)$. Now $g^{-1}\varphi (x)=x$,  
so $g^{-1}\varphi \in  \langle (y\,z), -\tau_x \rangle$
in view of Lemma \ref{lem:isostab}.  Thus,  $\varphi$ is
in the subgroup of  $\hbox{\rm Isom}_{\mathbb Z}(L)$ generated by
$G$, $(y\,z)$, $-\tau_x$.   But $\tau_x \in G$ so
$\varphi$ belongs to  the subgroup of $\hbox{\rm Isom}_{\mathbb Z}(L)$ generated by
$G$, $(y\,z)$, $-1$. Therefore $\hbox{\rm Isom}_{\mathbb Z}(L)$
is generated by $G$, $(y\,z)$, $-1$.   By the above comments,   
$\hbox{\rm Isom}_{\mathbb Z}(L)=
  \langle (y\,z), -1 \rangle \ltimes G$.
 The other assertions follow by symmetry or a routine argument. \qed
\bigskip

 \begin{theorem} \label{thm:equal}
 $\hbox{\rm AAut}_{\mathbb Z}(L)=\hbox{\rm Isom}_{\mathbb Z}(L)$.
\end{theorem}
\noindent {\it Proof:}   We know already  that 
$\hbox{\rm AAut}_{\mathbb Z}(L)\subseteq \hbox{\rm Isom}_{\mathbb Z}(L)$.
By Theorem 
\ref{thm:isostab} we have
$\hbox{\rm Isom}_{\mathbb Z}(L) =  \langle  (y\,z),-1 \rangle \ltimes G$.
But $(y\,z)$, $-1$, $G$  are all  contained in 
 $\hbox{\rm AAut}_{\mathbb Z}(L)$, so $ \hbox{\rm Isom}_{\mathbb Z}(L)
 \subseteq \hbox{\rm AAut}_{\mathbb Z}(L)$ holds as well.   \qed
 \bigskip

\begin{theorem}\label{thm:autL}
$\hbox{\rm Aut}_{\mathbb Z}(L)$ is equal to each of the groups
\begin{eqnarray*}
  \langle -(x\,y) \rangle \ltimes G,
  \qquad 
 \langle -(y\,z) \rangle \ltimes G,
  \qquad  \langle -(z\,x) \rangle \ltimes G.
\end{eqnarray*}
In particular $\hbox{\rm Aut}_{\mathbb Z}(L)$ is isomorphic to
${\mathbb Z}_2 \ltimes G$.
\end{theorem}
\noindent {\it Proof:} Combine 
\eqref{eq:normal}, Theorem
\ref{thm:isostab}, and Theorem \ref{thm:equal}. \qed
\medskip

\section {The lattice $L^* = \mathbb Z x^* \oplus  \mathbb Z y^*  \oplus \mathbb Z z^*$}  

In previous sections we have discussed the lattice
$L = \mathbb Z x \oplus  \mathbb Z y  \oplus \mathbb Z z$.
Here we consider the lattice
\begin{eqnarray}\label{eq:ls}
L^*: = \mathbb Z x^* \oplus  \mathbb Z y^*  \oplus \mathbb Z z^*
\end{eqnarray} 
from a similar point of view.
By 
\eqref{eq:shownil},  
$L^*$ is equal to $\mathfrak{sl}_2(\mathbb Z)$
and contains $L$.
Regarding $L$ and  $L^*$ as free abelian
groups,  we see {f}rom  \eqref{eq:mod2}  that $L$ has 
 index 4 in $L^*$  and
$L^*/L \cong \mathbb Z_2 \times  \mathbb Z_2$.
The duality between $\{x,y,z\}$ and  $\{x^*, y^*, z^*\}$ gives  \begin{eqnarray}
\label{eq:lsdef}
L^* &=& \lbrace u \in 
 \mathfrak{sl}_2
\,\vert \, (u,v) \in 2 \Z\quad \forall v \in L\rbrace,
\\
\label{eq:ldef}
L &=& \lbrace u \in 
 \mathfrak{sl}_2
\,\vert \, (u,v) \in 2 \Z\quad \forall v \in L^*\rbrace.
\end{eqnarray}
Define
\begin{eqnarray*}
\hbox{\rm Aut}_{\mathbb Z}(L^*) 
&=&
 \{ \varphi \in   \hbox{\rm Aut}_{\mathbb F}(\mathfrak {sl}_2) \mid   
\varphi(L^*) = L^*\},
\\
\hbox{\rm AAut}_{\mathbb Z}(L^*) 
&=&
 \{ \varphi \in   \hbox{\rm AAut}_{\mathbb F}(\mathfrak {sl}_2) \mid   
\varphi(L^*) = L^*\},
\\
\hbox{\rm Isom}_{\mathbb Z}(L^*) 
&=&
 \{ \varphi \in   \hbox{\rm Isom}_{\mathbb F}(\mathfrak {sl}_2) \mid   
\varphi(L^*) = L^*\}.
 \end{eqnarray*}   \smallskip

\begin{theorem}
\label{thm:starsame} We have
\begin{itemize}
\item[(i)]
$\hbox{\rm Aut}_{\mathbb Z}(L^*) 
=\hbox{\rm Aut}_{\mathbb Z}(L)$,
\item[(ii)]
$
\hbox{\rm AAut}_{\mathbb Z}(L^*) 
=\hbox{\rm AAut}_{\mathbb Z}(L),
$
\item[(iii)]
$\hbox{\rm Isom}_{\mathbb Z}(L^*) 
=\hbox{\rm Isom}_{\mathbb Z}(L)$.
\end{itemize}
\end{theorem}
\noindent {\it Proof:}  
By 
(\ref{eq:lsdef}) 
and (\ref{eq:ldef}), for each isometry $\sigma$ of 
$\mathfrak{sl}_2$ we have
\begin{eqnarray}
\sigma(L)=L \quad \Longleftrightarrow  \quad \sigma(L^*)=L^*.
\label{eq:iff}
\end{eqnarray}
Part (iii)
is immediate from this.
Parts (i) and (ii) also follow, since
by Lemma
\ref{eq:ginv} and
(\ref{eq:aainfo}) each of the groups
$\hbox{\rm Aut}_{\mathbb F}(\mathfrak {sl}_2)
$,
$\hbox{\rm AAut}_{\mathbb F}(\mathfrak {sl}_2)$
is contained in 
 $\hbox{\rm Isom}_{\mathbb F}(\mathfrak {sl}_2)$.
   \qed

\bigskip
Since  $\varrho \in G$ cyclically permutes $x^*, y^*,  z^*$, it follows that 
$G(x^*) = G(y^*) = G(z^*)$.   We will describe the orbit $G(z^*)$  
after first determining the stabilizer of $z^*$ in $G$. 
\bigskip

\begin{theorem} \label{thm:z*stab}  
The stabilizer of $z^*$ in $G$ is the subgroup generated
by $\exp(\ad z^*)$.  \end{theorem}

\noindent {\it Proof:} The subgroup of $G$
generated by $\exp(\ad z^*)$ is clearly contained in the stabilizer
of $z^*$.   

For the reverse containment, we observe by  \eqref{eq:thetahat} that to any $g \in G$ there corresponds 
$\theta \in \hbox{\rm SL}_2(\mathbb Z)$ such that ${\widehat \theta}=g$.
Writing 
$\theta= 
 \left [\begin{array}{cc}  m& p \\ n &  q \end{array} \right]
$
and using $z^* =
 \left [\begin{array}{cc}  0& -1 \\ 0 &  \ \ 0 \end{array} \right]
$,  we see that 
\begin{equation}\label{eq:gz*}
 g(z^*) = \theta z^* \theta^{-1}
= 
 \left [\begin{array}{cc}  mn& -m^2 \\ n^2 &  -mn \end{array} \right].
\end{equation}
When $g(z^*) = z^*$,  we must have  
$n=0$, which implies that $mq = \det(\theta) = 1$.    Hence $m=q \in  
\lbrace 1,-1\rbrace$.    Then  either $\theta$ or $-\theta$ is an integral power
of $\exp(z^*)$,   putting $g$ in the subgroup
generated by $\exp(\ad z^*)$. \qed 
\medskip

\begin{theorem}\label{thm:z*orbstab}  
\begin{itemize}   
\item[(i)]  The orbit $G(z^*)$  consists of $y^*, z^*$
together with all matrices of the form  $ \left[\begin{array}{cc}  mn & -m^2 \\
n^2 & -mn \end{array} \right]$
 where $m,n \in {\mathbb Z}$ are nonzero and relatively prime.
\item[(ii)] 
 $G(z^*)$ consists of $x^*, y^*, z^*$ together with all vectors of the form 
\begin{equation*}- \frac{1}{2}(a^2 x + b^2y + c^2 z),\end{equation*} 
where  $a,b,c$ are relatively prime positive integers
such that one  of  the relations  $a = b+c$, $b = c+a$, $c = a+b$ holds.
\end{itemize}
\end{theorem}  

\noindent {\it Proof:}  
For statement (i), we assume $g \in G$ and $g(z^*)$ is as in 
\eqref{eq:gz*} above.    Since  $\det(\theta) = mq-np=1$, 
either  $m,n$ are nonzero and relatively
prime,  or one of integers $m,n$ is zero.
If $m=0$ then $n=-p \in  
\lbrace 1,-1\rbrace$ and  $g(z^*)=y^*$.   We have
seen in the proof of Theorem \ref{thm:z*stab} that
when $n = 0$ then $g(z^*) = z^*$.   Hence (i) holds.  

Part (ii) expresses a matrix from part (i) 
in terms of the equitable basis.
This conversion to a combination of  $x,y,z$ just amounts to the identity
\begin{equation} \label{eq:orbelt}
\left[\begin{array}{cc}  mn & -m^2 \\
n^2 & -mn \end{array} \right]  = 
-\frac{1}{2} \Big ( (m-n)^2 x + m^2 y + n^2 z\Big). 
\end{equation}
\qed
 
 \begin{theorem}\label{thm:z*autoorb}  The sets $G(z^*)$ and $-G(z^*)$ are disjoint
 and the following sets coincide:
 \begin{itemize}
 \item[{\rm (i)}]   $G(z^*) \cup \big (- G(z^*)\big)$
 \item[{\rm (ii)}]  the orbit of $z^*$ under $\hbox{\rm Aut}_{\mathbb Z}(L^*)$
 \item[{\rm (iii)} ]  the orbit of $z^*$ under
 $\hbox{\rm AAut}_{\mathbb Z}(L^*)  =
  \hbox{\rm Isom}_{\mathbb Z}(L^*)$.
  \end{itemize}  \end{theorem}

\noindent {\it Proof:}   That  $G(z^*)$ and $-G(z^*)$ are disjoint can be seen from
Theorem \ref{thm:z*orbstab}.   The fact that the sets in (i)-(iii) are all equal is a
direct consequence of Theorems \ref{thm:isostab},   \ref{thm:equal},  \ref{thm:autL},  and  \ref{thm:starsame}.    \qed 
   \medskip

\section{Connections with a hyperbolic Kac-Moody Lie algebra}
In this section we make explicit the relationship between the equitable
basis for $\mathfrak{sl}_2$ 
and the Kac-Moody Lie algebra $\mathfrak{g}=\mathfrak{g}({\mathcal A})$ over $\mathbb F$
associated with the Cartan matrix
$\mathcal A$  from \eqref{eq:CM}. (All the terminology and necessary background material used here can be found in \cite[Ch.~1]{K}.)
The Coxeter-Dynkin diagram corresponding to $\mathcal A$ is
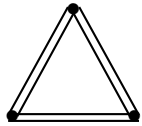
\begin{center}\begin{pspicture}(-2,-1)(1,1) 
\psset{xunit=.15cm,yunit=.15cm} 
{\psline{-}(-5,-3.2)(.2,6.3)}
{\psline{-}(-5.8,-3.1)(-.5,6.5)}
{\psline{-}(5,-3.2)(-.2,6.3)}
{\psline{-}(5.8,-3.1)(.5,6.5)}
{\psline{-}(-5.4,-3)(5.4,-3)}
{\psline{-}(-5.8,-3.6)(5.8,-3.6)}
\rput(-5.4,-3.2){{$\bullet$}}
\rput(5.4,-3.2){{$\bullet$}}
\rput(0,6.3){{$\bullet$}}    
%{\psline{-}(-2.3,4.6)(-1,5)}
%{\psline{-}(-1,5)(-.5,3.8)}
%{\psline{-}(2.3,4.6)(1,5)}
%{\psline{-}(1,5)(.5,3.8)}
%{\psline{-}(-4.6,-1.9)(-5.2,-.7)} 
%{\psline{-}(4.6,-1.9)(5.2,-.7)}  
%{\psline{-}(-4.6,-1.9)(-3.4,-1.4)}
%{\psline{-}(4.6,-1.9)(3.4,-1.4)}
%{\psline{-}(-3.5,-3.3)(-2.5,-2.3)}
%{\psline{-}(-3.5,-3.3)(-2.5,-4.3)}
%{\psline{-}(3.5,-3.3)(2.5,-2.3)}
%{\psline{-}(3.5,-3.3)(2.5,-4.3)}

\end{pspicture} 
\end{center} 
Since each subdiagram is  the diagram of a Cartan matrix of finite 
type $\hbox{\rm A}_1$ or of affine type $\hbox{\rm A}_1^{(1)}$, 
 the matrix $\mathcal A$ is of hyperbolic type (see \cite[\S 4.10]{K}).
We adopt the point of view that $x,y,z$ are the simple roots for $\mathfrak{g}$
and that 
there is a symmetric bilinear form $(\, ,\,)$ on $\hbox{\rm Span}_{\mathbb F}\lbrace x,y,z\rbrace$
whose values on these simple
roots are specified  by $\mathcal A$.  Since $\mathcal A$ is
nonsingular,  the simple roots are linearly independent  
and the bilinear form is nondegenerate.       
The root lattice 
may be identified with $L= \mathbb Z x \oplus \mathbb Z y \oplus \mathbb Z z$ 
 in our earlier notation.  Note that $L$ comes
equipped with the Lie product $\lbrack \ ,\ \rbrack$.
Note also that the group ${\rm Isom}_{\mathbb Z}(L)$ makes sense in
the present context. 
Three elements belonging to this group  
are the simple reflections  $r_x$, $r_y$, $r_z$.
For $u \in \lbrace x,y,z\rbrace$ we have $(u,u)=2$,   so the
reflection $r_u$ is given by 
\begin{equation}  r_u(v) = v - \frac{2(u,v)}{(u,u)} u =  v - (u,v) u \end{equation}
for all $v \in L$.   For example, 
\begin{eqnarray*}
r_x(x) = -x,  \qquad \qquad  r_x(y) = y + 2x, \qquad \qquad
r_x(z) = z + 2x. 
\end{eqnarray*}
Comparing this with Lemma \ref{lem3}\,(d) and using symmetry, we have
\begin{equation}\label{eq:reflects}
r_x = (y\,z)\tau_x, \qquad \qquad 
r_y = (z\,x)\tau_y, \qquad \qquad
r_z = (x\,y)\tau_z.
\end{equation} 
The subgroup $W$ of  ${\rm Isom}_{\mathbb Z}(L)$ generated by the
reflections  $r_x$, $r_y$, $r_z$  
is the {\it Weyl group}.
There is another subgroup of  ${\rm Isom}_{\mathbb Z}(L)$
that comes up naturally here.
Observe that $(x\,y)$, $(y\,z)$, $(z\,x)$ and $1$, $\varrho$, $\varrho^2$
together 
form a subgroup of
${\rm Isom}_{\mathbb Z}(L)$ that is isomorphic to the symmetric group
$S_3$; we identify this group with $S_3$ for the rest of the paper. 
In what follows, we adopt
$\pm S_3$ as a shorthand for $\langle \pm 1\rangle \times S_3$.   
We note that $S_3$ is the group of  diagram automorphisms
associated with $\mathcal A$ in the sense of
\cite[p.~68]{K}.
The following theorem is implied by
\cite[Cor.~5.10\,(b)]{K}.  \bigskip

\begin{theorem}
\label{thm:kacc}
 For the Cartan matrix
$\mathcal A$ in \eqref{eq:CM},
\begin{eqnarray*}
{\rm Isom}_{\mathbb Z}(L)= 
 \pm S_3 \ltimes W.
 \end{eqnarray*}
 \end{theorem}
 \bigskip

We now describe how $W$ is related to $G$.
Instead of doing this directly, we will relate
them both to a certain normal subgroup of
 $W$ denoted  $W^+$.
Recall that for $w \in W$ the {\it length} of $w$
is the number of factors 
 $r_x, r_y, r_z$ in
 a reduced expression for $w$. 
Let $W^+$ denote the subgroup of $W$ consisting of the
 elements of even
length.   Then  $W^+$ is a normal subgroup of $W$
with index 2. 

\begin{proposition}\label{prop:w+}
 $W^+=W \cap 
{\rm Aut}_{\mathbb Z}(L)$. 
Moreover $W^+ $ is a normal subgroup of 
${\rm Isom}_{\mathbb Z}(L)$ with index 24.
\end{proposition}
\noindent {\it Proof:}
To get the first assertion, note that each of
$r_x, r_y, r_z$ is an antiautomorphism
of $\mathfrak{sl}_2$ preserving $L$. The second assertion
follows from the first, using
Theorem \ref{thm:kacc} and the fact that
${\rm Aut}_{\mathbb Z}(L)$ is normal in
${\rm Isom}_{\mathbb Z}(L)$ with index 2.
 \qed \bigskip

The following result is immediate from the definition
of $W^+$.

 \begin{proposition}  
$W $ is equal to each of the groups
\begin{eqnarray*}
  \langle  r_x \rangle \ltimes W^+,
  \qquad 
 \langle r_y \rangle \ltimes W^+,
  \qquad  \langle r_z \rangle \ltimes W^+.
\end{eqnarray*}
In particular $W$ is isomorphic to 
 ${\mathbb Z}_2  \ltimes W^+$.
\end{proposition}  

\begin{proposition}
 $W^+$ is a normal subgroup of $G$, and 
the cosets of $W^+$ in $G$ are
\begin{eqnarray}
\label{eq:cosets}
W^+, \quad 
\tau_x W^+, \quad 
\tau_y W^+, \quad 
\tau_z W^+, \quad 
\varrho W^+, \quad 
\varrho^2 W^+.
\end{eqnarray}
The quotient group $G/W^+$ is isomorphic to $S_3$.
\end{proposition}
\noindent {\it Proof:}
The group $W^+$ is contained in $G$ since
\begin{equation}\label{eq:even}
 r_xr_y = \varrho \tau_z \tau_y,  \qquad  \quad 
r_y r_z = \varrho \tau_x \tau_z, \qquad
\quad r_z r_x = \varrho \tau_y \tau_x.
\end{equation}
Moreover, $W^+$ is normal in $G$ because $W^+$ is normal in 
${\rm Isom}_{\mathbb Z}(L)$.
The index of $W^+$ in $G$ is 6,  since
the index of $W^+$ in 
${\rm Isom}_{\mathbb Z}(L)$ is 24 and the index of
$G$ in 
${\rm Isom}_{\mathbb Z}(L)$ is 4.
Using 
(\ref{eq:reflects}) and 
(\ref{eq:even}),  one can readily verify
 that the list
(\ref{eq:cosets}) consists of the cosets of $W^+$ in $G$, and
that $G/W^+ \simeq S_3$.
\qed
\bigskip

\begin{proposition}\label{prop:isoWeyl}   $\hbox{\rm Isom}_{\mathbb Z}(L)/W^+  \cong \langle \pm 1 \rangle   \times D$, Πwhere $D$ is the dihedral
group of order 12.   \end{proposition}

\noindent {\it Proof:}  Let  $E = \hbox{\rm Isom}_{\mathbb Z}(L)/W^+$.   To argue that 
$E \cong \langle \pm 1 \rangle   \times D$,  we will produce elements $\eta, \vartheta \in E$ such
that $\eta$ has order 2, $\vartheta$ has order 6, and 

\begin{equation}\label{eq:1} \eta \vartheta \eta = \vartheta^{-1}.\end{equation}

\noindent  Set $\eta = (y \, z)W^+$ and $\vartheta = (y\,z)\tau_yW^+ = \tau_z(y \, z)W^+$.    
Note that $\eta$ has order 2, and since

$$(y \, z)\Big( (y \, z) \tau_y \Big)(y \, z) =   \tau_y (y\, z)   = \Big((y \, z) \tau_y \Big)^{-1},$$

\noindent equation  \eqref{eq:1} holds.     Moreover,

\begin{eqnarray*}  \Big ((y \, z)  \tau_y\Big)^2  &=&  (y \, z) \tau_y (y \, z)  \tau_y  \\
&=& \tau_z  \tau_y    \\
&=&(x \, y) r_z  (z\,x) r_y  \\
&=&(x \, y)(z\,x)  r_x  r_y \\
&=& \varrho^2 r_x r_y \\ 
& \equiv & \varrho^2   \quad  \mod W^+   \end{eqnarray*}

\noindent so that  $\vartheta^6 = (\varrho^2W^+)^3  = W^+$.   Thus,
$\vartheta$ has order 6.        Together $\eta, \vartheta$ generate a subgroup
of $E$ isomorphic to $D$.       Since
$\langle \pm 1\,W^+ \rangle$ is a central subgroup of $E$ and $|\,E\,| = 24$
by Proposition \ref{prop:w+}, we have
that $E \cong  \langle \pm 1 \rangle   \times D$, as claimed.     \qed

\bigskip

Let  $\Delta$ denote the set of roots attached to the Cartan matrix $\mathcal A$.
 Then $\Delta = \Delta_+ \cup \Delta_-$  where $\Delta_+$ (resp. $\Delta_-$) is the
 set of roots that are nonnegative (resp. nonpositive) integral linear combinations of
 the simple roots $x,y,z$. We have
$\Delta_-=-\Delta_+$. 
  A root is {\it real} if it lies in the  $W$-orbit of a simple root; 
otherwise it is {\it imaginary}.   Thus $\Delta$
decomposes into the disjoint union of the sets of real and imaginary roots:
$\Delta = \Delta^{\hbox{\rm re}} \cup
\Delta^{\hbox{\rm im}}$.  
The following theorem is implied by
\cite[Prop. 5.10\,(a)]{K}.

\bigskip
\begin{theorem}\label{thm3} For the Cartan matrix $\mathcal A$ 
in  \eqref{eq:CM} the set of real roots  $\Delta^{\hbox{\rm re}}$ coincides with the
set $R$ from Definition
\ref{def:r}.
\end{theorem}
\bigskip

With the above theorem in mind, we have  the following result  which gives an interpretation
of Proposition
\ref{prop4}.  \bigskip

\begin{proposition}\label{prop5} For the Cartan matrix $\mathcal A$
from \eqref{eq:CM},  the corresponding Weyl group $W$ and simple roots  satisfy
\begin{itemize}
\item[(i)]  $W(x) = \{ \alpha x + \beta y + \gamma z \in \Delta^{\hbox{\rm re}}  \mid
\alpha \equiv 1, \  \beta \equiv 0, \ \gamma \equiv 0  \mod 2\}$;
\item[(ii)]  $W(y) = \{ \alpha x + \beta y + \gamma z \in \Delta^{\hbox{\rm re}}  \mid
\beta \equiv 1, \  \alpha\equiv 0, \ \gamma \equiv 0  \mod 2\}$;
\item[(iii)]  $W(z) = \{ \alpha x + \beta y + \gamma z \in \Delta^{\hbox{\rm re}}  \mid
\gamma \equiv 1, \  \alpha\equiv 0, \  \beta \equiv 0  \mod 2\}$.
\end{itemize} 
\end{proposition}

\noindent {\it Proof:}   Each element in $W(x)$ is obtained from  $x$ by applying a product of reflections
in the simple roots.   The root $x$ belongs to the set on the right-hand side of (i), and whenever $u = 
\alpha x + \beta y + \gamma z$ belongs to that set,  then so do
\begin{eqnarray*}  r_x(u)  &=&  (2\beta + 2 \gamma -\alpha)x + \beta y + \gamma z, \\
r_y(u) &=&  \alpha x + (2 \alpha + 2 \gamma - \beta)y + \gamma z, \\
r_z(u) &=&  \alpha x  + \beta y + (2 \alpha +2 \beta - \gamma)z.  \end{eqnarray*} 
Consequently $W(x)$ is contained in the right side of (i).  Similar results apply in parts  (ii)
and (iii).  Thus  $\Delta^{\hbox{\rm re}}= W(x) \cup W(y) \cup W(z)$ is contained in the
(disjoint) union of the three sets on the right, which forces equality to hold in each case.   \qed

\bigskip

Consider the set of imaginary roots $\Delta^{\hbox{\rm im}}$ associated with
the Cartan matrix $\mathcal A$ from \eqref{eq:CM}.
It follows from \cite{M} or \cite[Prop.~5.2]{K} that
\begin{equation} \Delta^{\hbox{\rm im}} = \{ u \in \Delta \mid (u,u) \leq 0\}. \end{equation}
An important special case is the set of isotropic roots
\begin{equation} \Delta^0 = \{u \in \Delta \mid (u,u) = 0\}.\end{equation}
This set will be our focus for the rest of this section.

Each of the following four propositions contains a characterization 
$\Delta^0$.  \bigskip

\begin{proposition}
\label{prop:iso1pre}
 {\rm \cite[Prop. 5.10\,(c)]{K}}
The set of isotropic roots corresponding to the 
Cartan matrix $\mathcal A$ from \eqref{eq:CM} is given by
\begin{eqnarray*}
\Delta^0 = \lbrace u \in L\setminus \lbrace 0 \rbrace \;|\; (u,u)=0\rbrace.
\end{eqnarray*}
\end{proposition} \bigskip

The following result is implied by \cite[Prop.~5.7]{K}.  \bigskip

\begin{proposition}\label{prop:iso1} For the Cartan matrix $\mathcal A$
in \eqref{eq:CM},  a vector  is an isotropic root if and only if it is 
$W$-equivalent to a nonzero integer multiple of at least one of
\begin{eqnarray*}
 2x^* = -(y+z), \qquad  2y^* = -(z+x), \qquad 2z^* = -(x+y).
\end{eqnarray*}  
\end{proposition}    \bigskip

\begin{proposition}\label{prop:iso2}   For the Cartan matrix $\mathcal A$ in \eqref{eq:CM},
a vector is an isotropic root if and only if it is
$\hbox{\rm Isom}_{\mathbb Z}(L)$-equivalent to a positive even integer multiple of $z^*$
if and only if it is $G$-equivalent to a nonzero even integer multiple of $z^*$.   Thus, 
the set of isotropic roots is given by   

$$\Delta^0 = \bigcup_{n \in \mathbb Z, n \neq 0}    2nG (z^*).$$

\end{proposition}  
\noindent {\it Proof:}    This follows from Proposition \ref{prop:iso1},
{f}rom the fact that $\hbox{\rm Isom}_{\mathbb Z}(L) = \pm S_3 \ltimes W =  \langle (y\,z), -1 \rangle \ltimes G$,
and {f}rom the fact that 
$\Delta^0$ is $G$-invariant by Proposition \ref{prop:iso1pre}.  
  \qed   \bigskip

%Combining Theorem
%\ref{thm:z*orbstab}\,(iii) and Proposition
%\ref{prop:iso2} we routinely obtain the following result.

%\noindent {\it Proof:}    The first assertion follows readily from Proposition \ref{prop:iso1} and the fact %that that the generators  $\varrho$  and $\tau_y = (x\,z)r_y$  of $G$ leave $\Delta^0$ invariant. 
%The second sentence is just a restatement of the first using 
%Theorem \ref{thm1}.  
%  \qed   \bigskip
  
%  In Proposition \ref{prop:iso1} we saw that the isotropic roots for
%the Cartan matrix $\mathcal A$  are integer multiples of
%elements in $W(-2x^*) \cup W(-2y^*) \cup W(-2z^*)$.  
%Next we determine  information about the coefficients
%when isotropic roots are expressed in terms of the equitable
%basis.    \medskip

\begin{proposition}\label{prop:iso3} 
The isotropic roots corresponding to the Cartan matrix
$\mathcal A$ in \eqref{eq:CM} are precisely the vectors of the form
$n( a^2 x + b^2 y + c^2 z)$ for some $n \in \mathbb Z \setminus \{0\}$ and
 $a,b,c \in \mathbb Z_{\geq 0}$ (not all 0) 
such that at least one of 
\begin{eqnarray}
\label{eq:3pos}
a=b+c, \qquad \qquad 
b=c+a,\qquad \qquad
c=a+b.
\end{eqnarray}
    \end{proposition}   

%%%%%%%%%%%%save aug 10
%   An element $u \in \mathbb Z x + \mathbb Zy + \mathbb Zz$ is
%an isotropic root for the Cartan matrix $\mathcal A$ if only if  
%$u = n( \alpha^2 x + \beta^2 y + \gamma^2 z)$ for some $n \in \mathbb Z \setminus \{0\}$ and some %$\alpha,\beta,\gamma \in \mathbb Z_{\geq 0}$ not all 0 
%such that one of the following holds: \smallskip
%
%$\hbox{\rm (i)} \ \alpha = \beta + \gamma$,  \qquad \hbox{\rm (ii)}  $\beta = \gamma +\alpha$, \quad  or 
%\quad \hbox{\rm (iii)}  $\gamma = \alpha + \beta$.    \end{proposition}   
%%%save above aug10

\noindent {\it Proof:}  
Combine Theorem
\ref{thm:z*orbstab}\,(iii) and Proposition
\ref{prop:iso2}.
 \qed 
\bigskip

 Suppose now that $e_i,f_i,h_i \ (1 \leq i \leq 3)$ are the Chevalley generators
 for the Kac-Moody Lie algebra $\mathfrak g$ over 
$\mathbb F$
 corresponding to
 the Cartan matrix $\mathcal A$.  They satisfy the Serre relations (see \cite[(0.3.1)]{K}).      The Lie algebra has a decomposition 
 $\g = \mathfrak h \oplus \bigoplus_{u \in \Delta}  \mathfrak g_u$,  where
 $\mathfrak h$ is the span of the $h_i$.   For $u = 
 \alpha x + \beta y + \gamma z  \in \Delta_+$, 
$\mathfrak g_u$   is the subspace of $\mathfrak g$ consisting of all products of
$e_i$ such that the number of $e_1$'s appearing is $\alpha$, $e_2$'s is $\beta$
and $e_3$'s is $\gamma$.   A similar statement applies for $u \in \Delta_-$ with
$f_i$'s replacing the $e_i$'s.    Each subspace $\mathfrak g_u$ is finite-dimensional,
and its dimension is said to be the {\it multiplicity} of the root $u$.    There is an automorphism of
$\mathfrak g$ interchanging $e_i$ and $f_i$ and sending $h_i$ to $-h_i$.    
Thus, the multiplicities of $u$ and $-u$ are the same.  
Since

\begin{eqnarray*}
\dim \mathfrak g_{\sigma u} = \dim \mathfrak g_{u} 
\end{eqnarray*}
\noindent holds for all roots $u$ and all $\sigma \in W$,  it follows that
the multiplicity of each real root is 1.             For the isotropic roots,  we conclude the following.
\bigskip

\begin{corollary}\label{cor:isomult}  For the Cartan matrix
 $\mathcal A$ in \eqref{eq:CM},  each isotropic root has multiplicity 1. \end{corollary}
 
\noindent {\it Proof:}  Let $u \in \Delta^0$.   By   Proposition \ref{prop:iso1},   we may assume that
$u$ is a nonzero integer multiple of $2x^*$, $2y^*$, or $2z^*$.   Since the situation is
cyclically symmetric,  and since $u$ and $-u$ have the same multiplicity,  there is no loss in generality
in assuming $u = n(x+y)$ for $n \in \mathbb Z_{>0}$.    Thus,  elements in $\mathfrak g_u$ are
commutators with $n$ factors equal to $e_1$ and $n$ factors equal to $e_2$.  Since no $e_3,f_3,h_3$ are involved and since the relations governing $e_i,f_i,h_i$ for $i=1,2$ are the same as for
the affine Lie algebra $\mathfrak g(\hbox{\rm A}_1^{(1)})$   corresponding to the matrix
$\left[\begin{array}{cc} \ \ 2 & -2 \\ -2 & \ \ 2 \end{array}\right ],$ 
 the multiplicity of $u$ is the same as the
multiplicity of $n(x+y)$ in $\mathfrak g(\hbox{\rm A}_1^{(1)})$, which is known to be 1 (see \cite[Cor.~7.4]{K},  for example).
\qed  
\medskip

\section{$\mathfrak{sl}_2$-modules and the equitable basis}

Let $\mathcal V$  denote a finite-dimensional $\mathfrak{sl}_2$-module. 
Let $\varphi: \mathfrak{sl}_2 \rightarrow \mathfrak{gl}(\mathcal V)$
be the representation afforded 
 by $\mathcal V$, so that $\varphi(u)(v) = u.v$ for
all $u \in \mathfrak{sl}_2$ and $v \in \mathcal V$.  For all nilpotent $u \in \mathfrak{sl}_2$, 
the map $\varphi (u) :{\mathcal V} \to {\mathcal V}$
is nilpotent;   therefore the map 
$\exp\big(\varphi (u)\big) = \sum_{n=0}^\infty  \varphi (u)^n/n!$ is well-defined.
We note that $\exp\big(\varphi (u)\big)$ is invertible with inverse
$\exp\big(-\varphi (u)\big)$, and that 
\begin{equation}\label{eq:cong} 
 \exp\big(\varphi (u)\big) \varphi (t) \exp\big(\varphi (u)\big)^{-1} = \exp\big(\ad \varphi (u)\big)\big(\varphi(t)\big) \end{equation}
for all $t \in \mathfrak{sl}_2$.  
Using the dual basis $\{x^*,y^*,z^*\}$ for $\mathfrak{sl}_2$,  
we define the maps 
\begin{eqnarray*} \label{eq:PT}   P &=&  \exp\big(\varphi (x^*)\big)\exp\big(\varphi (y^*)\big),\nonumber \\
T_x &=& \exp\big(\varphi (y^*)\big)\exp\big(\varphi (z^*)\big) \exp\big(\varphi (y^*)\big), \nonumber \\
T_y &=& \exp\big(\varphi (z^*)\big)\exp\big(\varphi (x^*)\big) \exp\big(\varphi (z^*)\big),\nonumber  \\
T_z&=& \exp\big(\varphi (x^*)\big)\exp\big(\varphi (y^*)\big) \exp\big(\varphi (x^*)\big). \nonumber \end{eqnarray*}
If $\varphi$ is the 
adjoint representation,  then these are just the maps $\varrho$, $\tau_x$, $\tau_y$, $\tau_z$
from Section 3.   The next results say that analogues of the relations in
Lemmas \ref{lem1} and \ref{lem2} hold for arbitrary representations $\varphi$.   \bigskip

\begin{lemma} \label{lem7}  Let  $\varphi: \mathfrak{sl}_2 \rightarrow \mathfrak{gl}(\mathcal V)$  be a finite-dimensional $\mathfrak{sl}_2 $-representation.    Then 
the corresponding maps
$P, T_x,T_y, T_z$ satisfy {\rm (i)}--{\rm (v)} below.
\begin{itemize}
\item [(i)]  $P\varphi(x)P^{-1} = \varphi(y),  \ \ P\varphi(y) P^{-1} = \varphi(z), \ \ P\varphi(z)P^{-1} =\varphi(x)$;
\item [(ii)] $P\varphi(x^*)P^{-1} = \varphi(y^*),  \ \ P\varphi(y^*)P^{-1} = \varphi(z^*), \ \ P\varphi(z^*)P^{-1} =\varphi(x^*)$;
\item [(iii)] $P \exp\left(\varphi(x^*)\right) P^{-1} = \exp\left(\varphi(y^*)\right),  \ \ P \exp\left(\varphi(y^*)\right) P^{-1} = \exp\left(\varphi(z^*)\right), \\ P \exp\left(\varphi(z^*)\right) P^{-1} = \exp\left(\varphi(x^*)\right)$;
\item [(iv)]  $PT_xP^{-1} =T_y,  \ \ PT_yP^{-1} = T_z, \ \ PT_zP^{-1} =T_x$;
\item [(v)] $T_x \varphi(x) T^{-1}_x = -\varphi(x)$, \ \ 
$T_x \varphi(y)T^{-1}_x = 2 \varphi(x)+\varphi(z)$, \ \
$T_x \varphi(z)T^{-1}_x = 2 \varphi(x)+\varphi(y)$.
\end{itemize}
\end{lemma}

\noindent {\it Proof:}  Using \eqref{eq:cong}, we have
\begin{eqnarray*}
P \varphi (x)P^{-1} &=&
\exp\big(\ad \varphi(x^*)\big)\exp\big(\ad \varphi(y^*)\big)\big(\varphi(x)\big)
\\
&=& \varphi 
\big(\exp (\ad x^*)\exp (\ad y^*)(x)\big).
\end{eqnarray*}
Recall $\exp (\ad x^*)\exp (\ad y^*) = \varrho$ 
by Lemma
\ref{lem1}(ii)
and
$\varrho(x)=y$
by \eqref{eq:rho}
 so $P\varphi(x)P^{-1} = \varphi(y)$.

The remaining assertions can be
deduced from similar arguments.     \qed

\bigskip
\begin{lemma} \label{cor5}
Let  $\varphi: \mathfrak{sl}_2 \rightarrow \mathfrak{gl}(\mathcal V)$ be
 a finite-dimensional $\mathfrak{sl}_2 $-representation.    Then 
the corresponding maps
$P, T_x,T_y, T_z$ satisfy {\rm (i)}--{\rm (vi)} below.
\begin{itemize}
\item[(i)]  $P =  \exp\big(\varphi(x^*)\big)\exp\big(\varphi(y^*)\big)  =  \exp\big(\varphi(y^*)\big)\exp\big(\varphi(z^*)\big)$

 $=  \exp\big(\varphi(z^*)\big)\exp\big(\varphi(x^*)\big)$; 
\item[(ii)] $T_x = \exp\big(\varphi(y^*)\big)\exp\big(\varphi(z^*)\big)\exp\big(\varphi(y^*)\big)$

$\ = \exp\big(\varphi(z^*)\big)\exp\big(\varphi(y^*)\big) \exp\big(\varphi(z^*)\big)$;
\item[(iii)] $T_y =  \exp\big(\varphi(z^*)\big) \exp\big(\varphi(x^*)\big) \exp\big(\varphi(z^*)\big)$ 

$\ =  \exp\big(\varphi(x^*)\big) \exp\big(\varphi(z^*)\big) \exp\big(\varphi(x^*)\big)$;
\item[(iv)] $T_z = \exp\big(\varphi(x^*)\big) \exp\big(\varphi(y^*)\big) \exp\big(\varphi(x^*)\big)$

$\ =  \exp\big(\varphi(y^*)\big) \exp\big(\varphi(x^*)\big)  \exp\big(\varphi(y^*)\big)$;
\item[(v)]  $P^3 = T_x^2 = T_y^2 = T_z^2$;
\item[(vi)] Let $C$ denote the map in {\rm (v)}. 
Then $C$ commutes with $\varphi(u)$
for all $u \in \mathfrak {sl}_2$. 
  \end{itemize}
 \end{lemma}

\noindent {\it Proof:}   Parts (i)--(iv) follow easily from Lemma \ref{lem7}\,(iii),  so
consider part (v). 
To prove  that $P^3=T_x^2$, in the left-hand side of
$T_x T_x=T_x^2$ evaluate the first factor (resp. second 
factor) using the first (resp. second) equation in (ii) above,
and simplify the result using
$P=\exp\big(\varphi(y^*)\big)\exp\big(\varphi(z^*)\big)$.
The equations
$P^3=T_y^2$ and $P^3=T_z^2$ are obtained similarly.
Concerning (vi), note by Lemma 
\ref{lem7}(i) that $P^3$ commutes with each of $\varphi(x)$, $\varphi(y)$, $\varphi(z)$
and recall $x,y,z$ is a basis for $\mathfrak{sl}_2$. 
    \qed
\medskip

\begin{remark} \rm Let $\hbox{\rm B}_3$ be Artin's braid group given by
generators $s_1, s_2$ and the relation\,  $s_1 s_2 s_1 =
s_2 s_1 s_2$.    It  can be seen from 
Lemma \ref{cor5} that  each finite-dimensional $\mathfrak{sl}_2$-representation  $\varphi:
\mathfrak{sl}_2 \rightarrow \mathfrak{gl}(\mathcal V)$
determines a representation of $\hbox{\rm B}_3$ given by
$$s_1 \mapsto  \exp\big(\varphi(x^*)\big),  \qquad  \qquad s_2 \mapsto  \exp\big(\varphi(y^*)\big).$$
(Of course we could replace the pair $x,y$ by $y,z$ or $z,x$ in the above line.)

  The center of $\hbox{\rm B}_3$ is
generated by $(s_1 s_2)^3$ and this maps to $P^3 = T_x^2 = T_y^2 =
T_z^2$.        
In Corollary
\ref{cor:pttt} 
below,  we will show that $P^3$ must act as a scalar multiple of the identity  when $\mathcal V$ is
an irreducible
   $\mathfrak {sl}_2$-module, and we will determine the exact value of that scalar.
More general information on irreducible $\hbox{\rm B}_3$-modules, particularly
those of dimension $\leq 5$, can be found in \cite{TW}.  \end{remark} 

\medskip For each integer $d\geq 0$,
there exists
a unique irreducible $\mathfrak{sl}_2$-module of dimension
$d+1$ up to isomorphism.
This module, which we denote ${\mathcal V}(d)$,  has a  basis
$\lbrace v_i\rbrace_{i=0}^d$
% $v_0,v_1,\dots, v_d$ 
such that
% the  $\mathfrak{sl}_2$-action is
%given by  
\begin{eqnarray*} 
 h.v_i &=& (d-2i) v_i, \\
f.v_i &=& (i+1) v_{i+1}, \\
e.v_i &=& (d-i+1)v_{i-1}
 \end{eqnarray*}
 for $0 \leq i \leq d$,
where $v_{-1}=0$ and $v_{d+1}=0$.
We call $\lbrace v_i\rbrace_{i=0}^d$ a {\it standard basis} of $\mathcal V(d)$.  
%Adopting the convention that $v_{-1} = 0 = v_{d+1}$, we have that   
By \eqref{eqbas}
we have
%the equitable and dual bases on the standard basis is given by
\begin{eqnarray*}\label{eq:eqact}    x.v_i &=& (d-2i) v_i,   \\
y.v_i &=& 2(d-i+1) v_{i-1} + (2i-d)v_i,   \nonumber \\
z.v_i &=& (2i-d)v_i -2(i+1)v_{i+1}    \nonumber
 \end{eqnarray*}
 for $0 \leq i \leq d$, and by 
\eqref{eq:shownil} 
we have
\begin{eqnarray*}\label{eq:staract} 
 x^*.v_i &=& (i-d-1)v_{i-1}+(d-2i) v_i + (i+1)v_{i+1},  \\
y^*.v_i &=& (i+1) v_{i+1},   \nonumber \\
z^*.v_i &=& (i-d-1)v_{i-1}   \nonumber \end{eqnarray*}
for $0 \leq i \leq d$.
Let 
$\varphi_d$ denote the $\mathfrak{sl}_2$ representation afforded by ${\mathcal V}(d)$.
\bigskip

\begin{lemma}
\label{lem:expfacts}
With respect to a standard basis for ${\mathcal V}(d)$,
\begin{itemize}
\item[(i)] the matrix representing $\exp\big(\varphi_d(y^*)\big)$ is lower
triangular with $(i,j)$ entry ${i \choose j}$ for $0 \leq j\leq i\leq d$;
\item[(ii)] the matrix representing $\exp\big(-\varphi_d(y^*)\big)$ is lower
triangular with $(i,j)$ entry $(-1)^{i-j}{i \choose j}$ for $0 \leq j\leq i\leq d$;
\item[(iii)] the matrix representing $\exp\big(\varphi_d(z^*)\big)$ is upper
triangular with $(i,j)$ entry $(-1)^{j-i}{d-i \choose j-i}$ for $0 \leq i\leq j\leq d$;
\item[(iv)] the matrix representing $\exp\big(-\varphi_d(z^*)\big)$ is upper
triangular with $(i,j)$ entry ${d-i \choose j-i}$ for $0 \leq i\leq j\leq d$.
\end{itemize}
\end{lemma}
\noindent {\it Proof:} This is a routine calculation using the definition of the exponential
and the actions of $y^*$, $z^*$ on the standard basis. \qed
\medskip

\begin{lemma}
\label{lem:txaction}
For a standard basis $\lbrace v_i\rbrace_{i=0}^d$ of ${\mathcal V}(d)$,
\begin{eqnarray}
T_x v_i = (-1)^i v_{d-i} \qquad \qquad (0 \leq i \leq d).
\label{eq:txact}
\end{eqnarray}
\end{lemma}
\noindent {\it Proof:} By construction $v_i$ is an eigenvector for $x$ with
eigenvalue $d-2i$.    This combined with the first equation of Lemma
\ref{lem7}\,(v) implies  that $T_x v_i$ is an eigenvector for $x$ with eigenvalue $2i-d$.
Therefore there exists $\alpha_i \in {\mathbb F}$ such that
$T_x v_i=\alpha_i v_{d-i}$. We show $\alpha_i = (-1)^i$.
One readily checks that $\alpha_0=1$ using the definition of $T_x$
and the data in Lemma
\ref{lem:expfacts}. For $1 \leq i \leq d$,  we apply the second equation of
Lemma \ref{lem7}\,(v) to the vector $T_x v_i$; this yields $\alpha_i=-\alpha_{i-1}$ after
a brief calculation.
By the above comments $\alpha_i=(-1)^i$ for $0 \leq i \leq d$ and the result
follows. \qed
\medskip

\begin{corollary}
\label{cor:pttt}
For the maps $P$, $T_x$, $T_y$, $T_z$ corresponding to $\varphi_d$,
\begin{eqnarray}
P^3= T^2_x = T^2_y= T^2_z=(-1)^dI.
\label{eq:minusone}
\end{eqnarray}
\end{corollary}
\noindent {\it Proof:} Let $\lbrace v_i \rbrace_{i=0}^d$ denote a standard basis
for ${\mathcal V}(d)$. By Lemma
\ref{lem:txaction}, for $0 \leq i \leq d$ we see that 
$T_x v_i=(-1)^iv_{d-i}$ and
$T_x v_{d-i}=(-1)^{d-i}v_i $ so
$T_x^2 v_i=(-1)^d v_i$. Therefore $T_x^2= (-1)^dI$.
The result follows in view of Lemma
\ref{cor5}\,(v).
\qed

\bigskip
Starting with a standard basis $\lbrace v_i\rbrace_{i=0}^d$ of ${\mathcal V}(d)$ and using
the  map $P$ corresponding to $\varphi_d$, we obtain three different bases for
${\mathcal V}(d)$:
\begin{eqnarray}
\label{eq:3bases}
\lbrace v_i\rbrace_{i=0}^d,
\quad \qquad 
\lbrace Pv_i\rbrace_{i=0}^d,
\quad \qquad 
\lbrace P^2v_i\rbrace_{i=0}^d.
\end{eqnarray}
One significance of these bases is that for $0 \leq i \leq d$ the
vector $v_i$ (resp. $Pv_i$, resp. $P^2v_i$) is an eigenvector for
$x$ (resp. $y$, resp. $z$) with eigenvalue $d-2i$; this can be checked
using Lemma \ref{lem7}\,(i).
Our next goal is to describe how the three bases
(\ref{eq:3bases}) are related. To this end the following lemma will be useful. 
\bigskip

\begin{lemma}
For a standard basis $\lbrace v_i\rbrace_{i=0}^d$ of ${\mathcal V}(d)$
and for the map $P$ associated
with $\varphi_d$,
\begin{eqnarray}
\exp\big(\varphi_d(y^*)\big)v_i &=& (-1)^{d-i}P^2 v_{d-i},
\label{eq:p2one}
\\
\label{eq:p2two}
\exp\big(-\varphi_d(z^*)\big)v_i &=& (-1)^{d-i}Pv_{d-i}
\end{eqnarray}
for $0 \leq i \leq d$.
\end{lemma}
\noindent {\it Proof:}
Equation \eqref{eq:p2one} follows from 
\begin{eqnarray*}
P^{-2}\exp\big(\varphi_d(y^*)\big)v_i &=&
(-1)^d P\exp\big(\varphi_d(y^*)\big)v_i \\
&=&
(-1)^d \exp\big(\varphi_d(y^*)\big)
\exp\big(\varphi_d(z^*)\big) \exp\big(\varphi_d(y^*)\big)
v_i \\
&=&
(-1)^d T_xv_i \\
&=&
(-1)^{d-i}v_{d-i}
\end{eqnarray*}
and line
\eqref{eq:p2two} can be derived similarly. 
\qed  \medskip

\begin{corollary}
\label{cor:threesum}
For a standard basis $\lbrace v_i\rbrace_{i=0}^d$ of ${\mathcal V}(d)$
and for the map $P$ associated
with $\varphi_d$,
\begin{eqnarray}
Pv_0= \sum_{i=0}^d v_i,
 \qquad
P^2v_0= \sum_{i=0}^d Pv_i,
 \qquad
(-1)^d v_0= \sum_{i=0}^d P^2v_i.
\label{eq:threesum}
\end{eqnarray}
\end{corollary}
\noindent {\it Proof:} 
To derive the equation on the left in \eqref{eq:threesum}, 
 set $i=d$ in
\eqref{eq:p2two} and evaluate the result using
Lemma
\ref{lem:expfacts}(iv).
The other two relations in \eqref{eq:threesum}
can be obtained similarly  using $P^3=(-1)^dI$.
\qed  \medskip

 \begin{proposition}\label{prop3}  For a standard basis $\lbrace v_i\rbrace_{i=0}^d$ of
$\mathcal V(d)$ and for the map $P$ corresponding to $\varphi_d$ the following
{\rm (i)}--{\rm (iii)} hold for $0 \leq i \leq d$.    
\begin{itemize}
\item[(i)] The image of $\big(\varphi_d(z^*)\bigr)^{d-i}$ on ${\mathcal V}(d)$ is equal to
each of the subspaces
\begin{eqnarray*}
 \hbox{\rm Span}_{\mathbb F}\{v_0, \dots, v_i \}, \qquad
\hbox{\rm Span}_{\mathbb F}\{Pv_d,\dots, Pv_{d-i}\}.
\end{eqnarray*}
\item[(ii)] 
The image of $\bigl(\varphi_d(x^*)\bigr)^{d-i}$ on ${\mathcal V}(d)$ is equal to
each of the subspaces
\begin{eqnarray*}
 \hbox{\rm Span}_{\mathbb F}\{Pv_0, \dots, Pv_i \}, \qquad
 \hbox{\rm Span}_{\mathbb F}\{P^2v_d,\dots, P^2v_{d-i}\}.
\end{eqnarray*}
\item[(iii)] The image of $\bigl(\varphi_d(y^*)\bigr)^{d-i}$ on ${\mathcal V}(d)$ is equal to
each of the subspaces
\begin{eqnarray*}
 \hbox{\rm Span}_{\mathbb F}\{P^2v_0, \dots, P^2v_i \}, \qquad
\hbox{\rm Span}_{\mathbb F}\{v_d,\dots, v_{d-i}\}.
\end{eqnarray*}
\end{itemize}
\end{proposition}
\noindent {\it Proof:}
(i)  \ Recall $z^*.v_i=(i-d-1)v_{i-1}$ 
so the image of $\big(\varphi_d(z^*)\bigr)^{d-i}$ on ${\mathcal V}(d)$ 
is $\hbox{\rm Span}_{\mathbb F}\{v_0, \dots, v_i \}$.
Also $\hbox{\rm Span}_{\mathbb F}\{v_0, \dots, v_i \}=
\hbox{\rm Span}_{\mathbb F}\{Pv_d,\dots, Pv_{d-i}\}$ in view of
Lemma
\ref{lem:expfacts}(iv) 
and \eqref{eq:p2two}.

 (ii), (iii):     Apply $P$ and $P^2$ to the equations in part (i).
\qed  \medskip
 
We have described how the three bases 
\eqref{eq:3bases}
are related.   To visualize this description 
it is helpful to draw some diagrams.  In these pictures,  the following convention will be adopted.
Given bases
$\lbrace w_i\rbrace_{i=0}^d$ and $\lbrace w'_i\rbrace_{i=0}^d$ for $\mathcal V(d)$, the display
below will mean that
 $\hbox{\rm Span}_{\mathbb F}\{w_0,\dots, w_i\} 
=
 \hbox{\rm Span}_{\mathbb F}\{w_d', \dots, w_{d-i}'\}$ 
for $0 \leq i \leq d$:
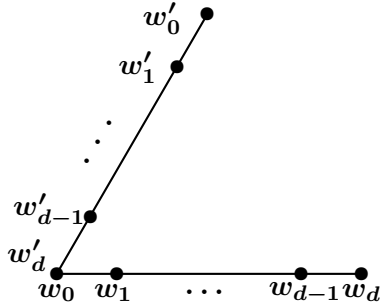
\medskip

 \begin{pspicture}(-4,-1)(4,4) 
\psset{xunit=.5cm,yunit=.5cm} 
{\psline{-}(0,0)(4,6.9)} 
{\psline{-}(0,0)(8,0)} 
\rput(8.1,-.5){${\black {\boldsymbol {w_d}}}$} 
\rput(6.6,-.5){${\black{\boldsymbol {w_{d-1}}}}$}   
\rput(4,-.5){${\black{\boldsymbol {\cdots}}}$}
\rput(1.5,-.5){${\black{\boldsymbol {w_1}}}$}
\rput(0,-.5){${\black{\boldsymbol {w_0}}}$} 
\rput(8.1,0){${\black {\boldsymbol {\bullet}}}$} 
\rput(6.5,0){${\black{\boldsymbol {\bullet}}}$}    
\rput(1.6,0){${\black{\boldsymbol {\bullet}}}$}
\rput(0,0){${\black{\boldsymbol {\bullet}}}$} 
 \rput(2.8,6.8){${\black {\boldsymbol {w_0'}}}$} 
\rput(2.15,5.5){${\black{\boldsymbol {w_{1}'}}}$}   
\rput(.8,3){${\black{\boldsymbol {\cdot} }}$}
\rput(1.1,3.5){${\black{\boldsymbol {\cdot }}}$}
\rput(1.4,4){${\black{\boldsymbol {\cdot} }}$}
\rput(-.2,1.7){${\black{\boldsymbol {w_{d-1}'}}}$}
\rput(-.8,.5){${\black{\boldsymbol {w_d'}}}$}  
 \rput(4,6.9){${\black {\boldsymbol {\bullet}}}$} 
\rput(3.2,5.5){${\black{\boldsymbol {\bullet}}}$}   
\rput(.9,1.5){${\black{\boldsymbol {\bullet}}}$} 
\end{pspicture}

\noindent With this convention in mind,
 Proposition \ref{prop3} tells us that

\bigskip
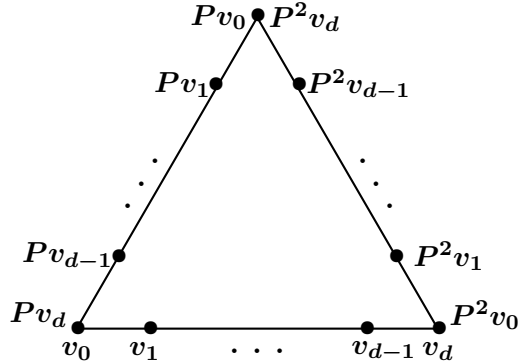
\begin{pspicture}(-4,-.5)(4,4.5) 
\psset{xunit=.6cm,yunit=.6cm} 
{\psline{-}(0,0)(4,6.9)}
{\psline{-}(4,6.9)(8,0)}
{\psline{-}(0,0)(8,0)} 
\rput(8,-.5){${\black {\boldsymbol {v_d}}}$} 
\rput(6.8,-.5){${\black{\boldsymbol {v_{d-1}}}}$}   
\rput(3.5,-.5){${\black{\boldsymbol {\cdot}}}$}
\rput(4,-.5){${\black{\boldsymbol {\cdot}}}$}
\rput(4.5,-.5){${\black{\boldsymbol {\cdot}}}$}
\rput(1.5,-.5){${\black{\boldsymbol {v_1}}}$}
\rput(0,-.5){${\black{\boldsymbol {v_0}}}$ 
\rput(7.7,.6){${\black {\boldsymbol {\bullet }}}$} 
\rput(6.1,.6){${\black{\boldsymbol {\bullet}}}$}   
\rput(1.3,.6){${\black{\boldsymbol {\bullet}}}$}
\rput(-.2,.6){${\black{\boldsymbol {\bullet}}}$ 
}
\rput(2.8,7.5){${\black {\boldsymbol {Pv_0}}}$} 
\rput(3.68,7.55){${\black {\boldsymbol {\bullet}}}$} 
\rput(2,6){${\black{\boldsymbol {Pv_1}}}$}   
\rput(2.75,6){${\black{\boldsymbol {\bullet}}}$}   
\rput(.8,3.3){${\black{\boldsymbol {\cdot} }}$}
\rput(1.1,3.8){${\black{\boldsymbol {\cdot }}}$}
\rput(1.4,4.3){${\black{\boldsymbol {\cdot} }}$}
\rput(-.55,2.2){${\black{\boldsymbol {Pv_{d-1}}}}$}
\rput(.6,2.2){${\black{\boldsymbol {\bullet}}}$}
\rput(-1.2,.9){${\black{\boldsymbol {Pv_d}}}$}
\rput(4.7,7.5){${\black {\boldsymbol {P^2v_d}}}$} 
\rput(5.9,6){${\black{\boldsymbol {P^2v_{d-1}}}}$}   
\rput(4.6,6){${\black{\boldsymbol {\bullet}}}$}   
\rput(6.6,3.3){${\black{\boldsymbol {\cdot} }}$}
\rput(6.3,3.8){${\black{\boldsymbol {\cdot} }}$}
\rput(6,4.3){${\black{\boldsymbol {\cdot}}}$}
\rput(7.9,2.2){${\black{\boldsymbol {P^2v_1}}}$}
\rput(6.75,2.2){${\black{\boldsymbol {\bullet}}}$}
\rput(8.7,.9){${\black{\boldsymbol {P^2v_0}}}$}}
\end{pspicture}

\medskip
\noindent By Corollary \ref{cor:threesum},  the sum of
the vectors on each edge of the triangle is a scalar multiple of the vector at
the opposite vertex. 
\medskip

For the special case of the adjoint module $\mathcal V(2)$,   the elements
\begin{eqnarray*}
v_0 = z^* = -e,  \qquad v_1 = x = h,  \qquad v_2 = y^* = f
\end{eqnarray*}
form a standard basis. In this case the picture is

\bigskip 
\begin{pspicture}(-4,0)(4,4) 
\psset{xunit=.5cm,yunit=.5cm} 
{\psline{-}(0,0)(4,6.9)}
{\psline{-}(4,6.9)(8,0)}
{\psline{-}(0,0)(8,0)} 
\rput(8.6,.1){${\black {\boldsymbol {y^*}}}$}  
\rput(4,-.6){${\black{\boldsymbol {x}}}$}
 \rput(-.7,.1){${\black{\boldsymbol {z^*}}}$ 
\rput(8.2,.25){${\black {\boldsymbol {\bullet }}}$}  
\rput(4.3,.25){${\black {\boldsymbol {\bullet }}}$}  
\rput(.35,.25){${\black{\boldsymbol {\bullet}}}$}  
\rput(1.7,3.9 ){${\black{\boldsymbol {y} }}$}  
\rput(6.8,3.9 ){${\black{\boldsymbol {z} }}$}  
\rput(4.5,7.8 ){${\black{\boldsymbol {x^*} }}$}  
\rput(2.4,3.8){${\black{\boldsymbol {\bullet} }}$}
\rput(4.3,7.2){${\black{\boldsymbol {\bullet} }}$} 
\rput(6.2,3.8){${\black{\boldsymbol {\bullet} }}$}
}
\end{pspicture}  \medskip

\bigskip    
\noindent By \eqref{eq:3sum} we have
\begin{eqnarray*}
  z^* + x + y^* = x^*,  \qquad 
x^* + y + z^* = y^*, \qquad y^* + z + x^* = z^*,
 \end{eqnarray*}  
 which is just \eqref{eq:threesum} in this special case.

\medskip
It is well known that each
irreducible $\mathfrak{sl}_2$-module $\mathcal V(d)$
can be realized explicitly as the space of homogeneous polynomials over
$\mathbb F$  of total degree
$d$  in two variables.     The symmetry in the equitable basis for $\mathfrak{sl}_2$ 
is reflected in how it acts on these polynomials, as we now discuss.  
\medskip

The action of $\mathfrak{sl}_2$ on $\mathcal V(1)$ extends to
an action of $\mathfrak{sl}_2$  on the symmetric algebra $\mathcal S: = S(\mathcal V(1))$ 
by derivations, so that $d.(uv) = (d.u)v + u(d.v)$ holds for all $d \in \mathfrak{sl}_2$
and all $u,v \in \mathcal S$.  We regard $\mathcal S$ as the polynomial algebra $\mathbb F[s,t]$
in two commuting indeterminates $s,t$
 and set $r = -s-t$.  Then 
$$\mathcal S = \mathbb F[r,s] = \mathbb F[s,t] = \mathbb F[t,r] \quad \hbox{and}
 \quad r+s+t = 0.$$  We identify the variables  $t,s$ with the standard basis elements $v_0,v_1$
 respectively, and get the actions 
 
\begin{eqnarray*}\label{eq:sact} \begin{array}{ccccccccc}
 x:&r\mapsto& s-t \quad  \quad &y:&s\mapsto& t-r \quad \quad   &z:&t\mapsto& r-s \\
    &s\mapsto& -s \quad  \quad  &&t\mapsto& -t \quad \quad     &&r\mapsto& -r \\
  &t\mapsto& t  \quad\quad  &&r\mapsto& r\quad \quad  &&s \mapsto& s, \end{array}
  \end{eqnarray*}
  \begin{eqnarray*}\label{eq:sact2} \begin{array}{ccccccccc}
 x^*: &r\mapsto&\ 0 \quad  \qquad  &y^*: &s\mapsto&\ 0\quad \qquad&z^*:&t\mapsto&\ 0 \\
    &s\mapsto&\ r \quad  \qquad   &&t\mapsto&\ s \quad \qquad&&r \mapsto&\ t \\
  &t\mapsto&-r \quad\qquad   &&r\mapsto& -s \quad \qquad&&s\mapsto& -t. \end{array}
  \end{eqnarray*}
  \medskip

\begin{proposition}\label{prop:symact}  For an integer $d \geq 0$ the following hold.
\begin{itemize}
\item [(i)]
 The
homogeneous polynomials of degree $d$ form an $\mathfrak{sl}_2$-submodule
  of  $\mathcal S$  that is isomorphic to $\mathcal V(d)$;

\item [(ii)] 
each of the sets $\lbrace s^i t^{d-i}\rbrace_{i=0}^d$,
 $\lbrace t^i r^{d-i}\rbrace_{i=0}^d$, $\lbrace r^i s^{d-i}\rbrace_{i=0}^d$
 is a basis for this submodule;
 
\item [(iii)]
 for $0 \leq i \leq d$ the vector 
$s^{i}t^{d-i}$ (resp.
$t^{i}r^{d-i}$, resp.  $r^{i}s^{d-i}$)
is an eigenvector for $x$ (resp. $y$, resp. $z$)
 with eigenvalue $d-2i$.
 \end{itemize}
 \end{proposition}

For $d=3$ these three bases appear on the perimeter of the triangle
below.
 
 \vskip .7 truein
  
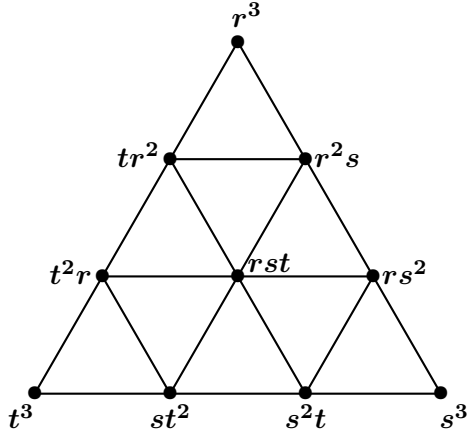
 \begin{pspicture}(-4,0)(4,4) 
\psset{xunit=.6cm,yunit=.6cm} 
{\psline{-}(0,0)(4.5,7.8)}
{\psline{-}(4.5,7.8)(9,0)}
{\psline{-}(0,0)(9,0)} 
\rput(9.3,-.5){${\black {\boldsymbol {s^3}}}$}  
\rput(6,-.5){${\black {\boldsymbol {s^2t}}}$}  
\rput(3,-.5){${\black{\boldsymbol {st^2}}}$}
\rput(-.3,-.5){${\black{\boldsymbol {t^3}}}$}

\rput(9,0){${\black {\boldsymbol {\bullet }}}$}  
\rput(6,0){${\black {\boldsymbol {\bullet }}}$}  
\rput(3,0){${\black {\boldsymbol {\bullet }}}$}  
\rput(0,0){${\black{\boldsymbol {\bullet }}}$} 
 
\rput(.8,2.7 ){${\black{\boldsymbol {t^2r}}}$}  
\rput(1.5,2.6){${\black{\boldsymbol {\bullet}}}$}
\rput(3,5.2){${\black{\boldsymbol {\bullet} }}$} 
\rput(2.3,5.3 ){${\black{\boldsymbol {tr^2}}}$}  
\rput(4.5,7.8){${\black{\boldsymbol {\bullet} }}$} 
\rput(4.7,8.4 ){${\black{\boldsymbol {r^3}}}$}  
\rput(6,5.2){${\black{\boldsymbol {\bullet} }}$}

\rput(6.7,5.3 ){${\black{\boldsymbol {r^2s}}}$}   
 \rput(7.5,2.6){${\black{\boldsymbol {\bullet} }}$}
 \rput(8.2,2.7 ){${\black{\boldsymbol {rs^2 }}}$}  
 {\psline{-}(1.5,2.6)(7.5,2.6)} 
{\psline{-}(1.5,2.6)(3,0)} 
{\psline{-}(3,0)(6,5.2)} 
 {\psline{-}(6,0)(3,5.2)} 
  {\psline{-}(6,0)(7.5,2.6)} 
 {\psline{-}(6,5.2)(3,5.2)} 
\rput(4.5,2.6){${\black{\boldsymbol {\bullet }}}$} 
\rput(5.2,2.9){${\black{\boldsymbol {rst }}}$} 
\end{pspicture}
\bigskip \bigskip  \medskip

  \section{Connections with the Poincar\'e Disk and Pythagorean Triples}
  
  Recall that $(v,w^*) = 2\delta_{v,w}$ for  $v,w \in \{x,y,z\}$.   Thus, the
  elements $\half x^*$, $\half y^*$, $\half z^*$ are the {\it fundamental weights}
  of the Kac-Moody algebra $\mathfrak g$,  and the lattice $\half L^* = \mathbb Z(\half x^*) \oplus
  \mathbb Z(\half y^*) \oplus \mathbb Z(\half z^*)$ is the {\it weight lattice}.  
  The elements in the set $\Omega := W(x^*) \cup W(y^*) \cup W(z^*)$
  are the Weyl group images of twice the fundamental weights.    Here we consider the
  set $\Omega$, and relate it to Pythagorean triples.     By  
   Proposition \ref{prop:iso1},
  Proposition \ref{prop:iso3}, and Theorem \ref{thm:z*orbstab},   the elements $u  \in \Omega$ have the form
  $u = -\frac{1}{2}\big( a^2 x + b^2 y + c^2 z\big)$
  where $a,b,c \in \mathbb Z_{\geq 0}$ (not all 0) and at least one of  the following holds:
  \begin{eqnarray*}
  a=b+c, \qquad \qquad b=c+a, \qquad \qquad c=a+b.
  \end{eqnarray*}
  Each element in $\Omega$ corresponds to a Pythagorean triple  $(\alpha,\beta, \gamma)$ in the following way.
  When $c = a+b$,  the Pythagorean triple $(\alpha, \beta,\gamma)$ can be
  obtained  from the matrix equation $M(a^2,b^2,c^2)^{\mathfrak t} = (\alpha,\beta,\gamma)^{\mathfrak t}$
  where 
  \begin{equation}\label{eq:changer} M = \left[\begin{array}{ccc}  -1 & \ \  1 &  1 \\ \ \  0 &  -1 & 1 \\ \ \ 0 & \ \ 1 & 1 \end{array} \right ], \end{equation}
  and $\gamma^2 = \alpha^2 + \beta^2$.    Similarly, when  $a = b+c$, then 
  for $(\beta,\gamma,\alpha)^{\mathfrak t}:= M(b^2,c^2,a^2)^{\mathfrak t}$ we have 
  $\alpha^2 = \beta^2 + \gamma^2$, and
  when    $b = c+a$ then for 
 $(\gamma, \alpha, \beta)^{\mathfrak t}: = M(c^2,a^2,b^2)^{\mathfrak t}$ we have
  $\beta^2 = \gamma^2 + \alpha^2$.

  There is a beautiful way to visualize the set $\Omega$  that brings together
  many of the ideas in this paper.  
  Starting with the picture for $\mathcal V(2)$ in Section 8, 
  and applying reflections in $W$ we obtain the Poincar\'e disk $\mathcal P$.  
  We have displayed a portion of the disk in Figure 1 below
  and have given the corresponding  triple $a^2, b^2, c^2$ for each displayed point on
  the circumference between $x^*$ and $y^*$.     For each of these points, $c = a +b$. 
  The other labels can be obtained by permuting $x^*,y^*,z^*$ and $x,y,z$.
  The vertices on the perimeter are the elements of $\Omega$;  that is, they are  the Weyl group reflections
  of the fundamental weights of the hyperbolic Kac-Moody algebra, after each
  one is multiplied by a factor of $2$.      
  
  The group   ${\rm Isom}_{\mathbb Z}(L) = 
\langle (y\, z), -1 \rangle \ltimes G  =  \pm S_3 \ltimes W$ is the
group of automorphisms  $\hbox{Aut}(\mathcal P)$ of the disk.   Indeed, the Weyl group
$W$ permutes the triangles.  By multiplying any $\xi  \in \hbox{Aut}(\mathcal P)$   
 by an appropriate element of $W$, we can assume 
that $\xi$  maps the central triangle  to itself.
Such an automorphism can be seen to belong to $\pm S_3$.
For the other realization of $ \hbox{Aut}(\mathcal P)$, observe that by multiplying 
an automorphism  $\xi  \in \hbox{Aut}(\mathcal P)$  by a suitable 
element of $G$, we can assume that $\xi$ fixes the edge labeled by
$x$ and $-x$.     Since the stabilizer of
$x$ in $G$ is trivial according to Lemma \ref{lem4},
such an automorphism must belong
to $\langle (y\, z), -1 \rangle$.

  \vfill \eject 
{\hskip 1.4 truein {\Large Figure 1}}
\vskip .5 truein 

\begin{pspicture}(-4.5,-6)(10,10) 
\psset{xunit=1cm,yunit=1cm} 
{\cnode(0,0){8}{b}}  
\rput(0,8){{$\bullet$}}
\rput(8,0){{$\bullet$}}
\rput(0,-8){{$\bullet$}}
\rput(-8,0){{$\bullet$}}
\rput(7.75,2){{$\bullet$}}
\rput(-7.75,2){{$\bullet$}}
\rput(7.75,-2){{$\bullet$}}
\rput(-7.75,-2){{$\bullet$}}
\rput(6.9,4){{$\bullet$}}
\rput(-6.9,4){{$\bullet$}}
\rput(6.9,-4){{$\bullet$}}
\rput(-6.9,-4){{$\bullet$}}
\rput(5.6,5.6){{$\bullet$}}
\rput(-5.6,5.6){{$\bullet$}}
\rput(5.6,-5.6){{$\bullet$}}
\rput(-5.6,-5.6){{$\bullet$}}
\rput(4,6.9){{$\bullet$}}
\rput(-4,6.9){{$\bullet$}}
\rput(4,-6.9){{$\bullet$}}
\rput(-4,-6.9){{$\bullet$}}
\rput(2,7.75){{$\bullet$}}
\rput(2,-7.75){{$\bullet$}}
\rput(-2,7.75){{$\bullet$}}
\rput(-2,-7.75){{$\bullet$}}
{\pscurve{-}(6.9,-4)(2.4,1.4)(0,8)}
{\pscurve{-}(-6.9,-4)(-2.4,1.4)(0,8)}
{\pscurve{-}(-6.9,-4)(0,-2.8)(6.9,-4)}
 \rput(2.1,1.2){\Large $\boldsymbol{z}$}
\rput(-2.1,1.2){\Large $\boldsymbol{y}$}
\rput(0,-2.3){\Large $\boldsymbol{x}$}
%{\pscurve{->}(-.2,-.4)(-.5,0)(0,.5)(.5,0)} 

 {\pscurve{-}(6.9,4)(2.3,5.1)(0,8)}
{\pscurve{-}(6.9,4)(5.5,0)(6.9,-4)}
\rput (2.7,1.6){\Large {$\boldsymbol{-z}$}} 
%{\pscurve{<-}(4.5,2.3)(4.1,2.6)(3.8,2.2)(4.1,1.8)}
 
{\pscurve{-}(-6.9,4)(-2.3,5.1)(0,8)}
{\pscurve{-}(-6.9,4)(-5.5,0)(-6.9,-4)}
\rput (-2.8,1.6){\Large {$\boldsymbol{-y}$}} 
%{\pscurve{<-}(-3.5,2.5)(-3.8,2.8)(-4.2,2.4)(-3.8,2)}
{\pscurve{-}(-6.9,-4)(-2.3,-5.1)(0,-8)}
{\pscurve{-}(6.9,-4)(2.3,-5.1)(0,-8)}
%{\pscurve{<-}(.3,-4.4)(-.1,-4.1)(-.4,-4.5)(-.1,-4.9)}
\rput(-.1,-3.2){\Large {$\boldsymbol{-x}$}}  
{\pscurve{-}(-6.9,4)(-5,5)(-4,6.9)}
{\pscurve{-}(6.9,4)(5,5)(4,6.9)}
{\pscurve{-}(-6.9,-4)(-5,-5)(-4,-6.9)}
{\pscurve{-}(6.9,-4)(5,-5)(4,-6.9)}
{\pscurve{-}(0,-8)(2.1,-6.7)(4,-6.9)}
{\pscurve{-}(0,-8)(-2.1,-6.7)(-4,-6.9)}
{\pscurve{-}(0,8)(-2.1,6.7)(-4,6.9)}
{\pscurve{-}(0,8)(2.1,6.7)(4,6.9)}
{\pscurve{-}(8,0)(6.7,2.1)(6.9,4)}
{\pscurve{-}(8,0)(6.7,-2.1)(6.9,-4)}
{\pscurve{-}(-8,0)(-6.7,2.1)(-6.9,4)}
{\pscurve{-}(-8,0)(-6.7,-2.1)(-6.9,-4)} 
 
 {\pscurve{-}(8,0)(7.5,1)(7.75,2)}
{\pscurve{-}(-8,0)(-7.5,1)(-7.75,2)}
{\pscurve{-}(-8,0)(-7.5,-1)(-7.75,-2)}
{\pscurve{-}(8,0)(7.5,-1)(7.75,-2)}
{\pscurve{-}(7.75,2)(6.95,3)(6.9,4)}
{\pscurve{-}(-7.75,2)(-6.95,3)(-6.9,4)}
{\pscurve{-}(-7.75,-2)(-6.95,-3)(-6.9,-4)}
{\pscurve{-}(7.75,-2)(6.95,-3)(6.9,-4)}
{\pscurve{-}(6.9,4)(5.8,4.8)(5.6,5.6)}
{\pscurve{-}(6.9,-4)(5.8,-4.8)(5.6,-5.6)}
{\pscurve{-}(-6.9,4)(-5.8,4.8)(-5.6,5.6)}
{\pscurve{-}(-6.9,-4)(-5.8,-4.8)(-5.6,-5.6)}
{\pscurve{-}(5.6,5.6)(4.8,5.8)(4,6.9)}
{\pscurve{-}(5.6,-5.6)(4.8,-5.8)(4,-6.9)}
{\pscurve{-}(-5.6,5.6)(-4.8,5.8)(-4,6.9)}
{\pscurve{-}(-5.6,-5.6)(-4.8,-5.8)(-4,-6.9)}
{\pscurve{-}(4,6.9)(3,6.95) (2,7.75)}
{\pscurve{-}(4,-6.9)(3,-6.95) (2,-7.75)}
{\pscurve{-}(-4,6.9)(-3,6.95) (-2,7.75)}
{\pscurve{-}(-4,-6.9)(-3,-6.95) (-2,-7.75)}
{\pscurve{-}(2,7.75)(1,7.5)(0,8)}
{\pscurve{-}(2,-7.75)(1,-7.5)(0,-8)}
{\pscurve{-}(-2,7.75)(-1,7.5)(0,8)}
{\pscurve{-}(-2,-7.75)(-1,-7.5)(0,-8)}
%\rput(-3.6,-2.3){\large{$\boldsymbol {\exp(\hbox{\rm \textbf {ad}}\,{\blue z^*})}$}}
%{\psset{arrows=->, linecolor=red}{{\pscurve{->}(-3.2,-1.2)(-2,-1.6)(-1.6,-2.5)} }}
%\rput(3.6,-2.3){\large{$\boldsymbol {\exp(\hbox{\rm \textbf {ad}}\,{\blue y^*})}$}}
%{\psset{arrows=->, linecolor=red}{{\pscurve{<-}(3.2,-1.2)(2,-1.6)(1.6,-2.5)} }}
%\rput(0,3.5){\large{$\boldsymbol {\exp(\hbox{\rm \textbf {ad}}\,{\blue x^*})}$}}
%{\psset{arrows=->, linecolor=red}{{\pscurve{->}(.9,3)(0,2.8)(-.9,3)} }} 
 
\rput(0.1,8.4){\large{$\boldsymbol{\blue x^*}$}}
\rput (0,9.4){\Large{$\boldsymbol{\red 0^2}$}}
\rput(0,10) {\Large{$\boldsymbol{\red 1^2}$}}
\rput(0,10.6) {\Large{$\boldsymbol{\red 1^2}$}}

 \rput(2.5,9.1){\Large{$\boldsymbol{\red 1^2}$}} 
  \rput(2.7,9.7){\Large{$\boldsymbol{\red 3^2}$}} 
   \rput(2.9,10.3){\Large{$\boldsymbol{\red 4^2}$}} 
   
  \rput(3.3,7.6){\large{$\boldsymbol{\blue 6}$}}
   \rput(3.65,7.55){\large{$\boldsymbol{\blue x^*}$}}
  \rput(3.9,7.35){\large{$\boldsymbol{\blue +}$}}
  \rput(4.2,7.15){\large{$\boldsymbol{\blue 3}$}}
   \rput(4.5,7.05){\large{$\boldsymbol{\blue y^*}$}}
    \rput(4.7,6.8){\large{$\boldsymbol{\blue -}$}}
     \rput(5.0,6.62){\large{$\boldsymbol{\blue 2}$}}
      \rput(5.3,6.6){\large{$\boldsymbol{\blue z^*}$}}

  \rput(5.1,8.3) {\Large{$\boldsymbol{\red 1^2}$}}
   \rput(5.5,8.9) {\Large{$\boldsymbol{\red 2^2}$}}
    \rput(5.9,9.5) {\Large{$\boldsymbol{\red 3^2}$}}

    \rput(6.8,6.8){\Large{$\boldsymbol{\red 2^2}$}}
      \rput(7.3,7.4){\Large{$\boldsymbol{\red 3^2}$}}
        \rput(7.8,8){\Large{$\boldsymbol{\red 5^2}$}}

    \rput(6.55,5){\large{$\boldsymbol{\blue 2}$}}
   \rput(6.9,4.95){\large{$\boldsymbol{\blue x^*}$}}
  \rput(7,4.6){\large{$\boldsymbol{\blue +}$}}
  \rput(7.2,4.2){\large{$\boldsymbol{\blue 2}$}}
   \rput(7.5,4.15){\large{$\boldsymbol{\blue y^*}$}}
    \rput(7.5,3.75){\large{$\boldsymbol{\blue -}$}}
     \rput(7.9,3.7){\large{$\boldsymbol{\blue z^*}$}} 
     
        \rput(-7.5,3.5){\large{$\boldsymbol{\blue z^*}$}} 
         \rput(-7.35,3.9){\large{$\boldsymbol{\blue +}$}}
       \rput(-7.8,3.5){\large{$\boldsymbol{\blue 2}$}}
   \rput(-7.1,4.3){\large{$\boldsymbol{\blue x^*}$}}  
     \rput(-7.4,4.3){\large{$\boldsymbol{\blue 2}$}}
     \rput(-7,4.75){\large{$\boldsymbol{\blue -}$}}
     \rput(-6.55,4.8){\large{$\boldsymbol{\blue y^*}$}}

\rput(8.4,5){\Large{$\boldsymbol{\red 1^2}$}}
\rput(9,5.5){\Large{$\boldsymbol{\red 1^2}$}} 
\rput(9.7,6){\Large{$\boldsymbol{\red 2^2}$}}

 \rput(9.1,2.6){\Large{$\boldsymbol{\red 3^2}$}}
  \rput(9.8,3){\Large{$\boldsymbol{\red 2^2}$}}
   \rput(10.5,3.4){\Large{$\boldsymbol{\red 5^2}$}}
   
 \rput(8.7,.25){\large{$\boldsymbol{\blue 3x^*}$}}
 \rput(9.6,.25){\large{$\boldsymbol{\blue+6y^*}$}}
 \rput(10.6,.25){\large{$\boldsymbol{\blue-2z^*}$}} 
 \rput (9.3,-.25){\Large{$\boldsymbol{\red 2^2}$}}
  \rput (10.1,-.25){\Large{$\boldsymbol{\red 1^2}$}}
   \rput (10.8,-.25){\Large{$\boldsymbol{\red 3^2}$}}

  \rput(9.1,-2.2){\Large{$\boldsymbol{\red 3^2}$}} 
  \rput(9.9,-2.4){\Large{$\boldsymbol{\red 1^2}$}}
\rput(10.7,-2.6){\Large{$\boldsymbol{\red 4^2}$}} 

  \rput(7.4,-4.2){\large{$\boldsymbol{\blue y^*}$}} 
    \rput(-7.4,-4.2){\large{$\boldsymbol{\blue z^*}$}}

  \rput(8.2,-4.6) {\Large{$\boldsymbol{\red 1^2}$}}
   \rput(8.9,-5) {\Large{$\boldsymbol{\red 0^2}$}}
    \rput(9.6,-5.5) {\Large{$\boldsymbol{\red 1^2}$}}

    \rput(-.9,-8.3){\large{$\boldsymbol{\blue 2y^*}$}}
 \rput(0,-8.3){\large{$\boldsymbol{\blue+2z^*}$}}
 \rput(.9,-8.3){\large{$\boldsymbol{\blue-x^*}$}}
 
\end{pspicture}

\vfill \eject


\begin{thebibliography}{ITW}
\bibitem[A1]{A1}  R.C.~Alperin,   $\hbox{\rm PSL}_2(\mathbb Z) = \mathbb Z_2 \ast \mathbb Z_3$,  Amer. Math. Monthly
\textbf{100} (1993),  385-386.  

\bibitem[A2]{A2}  R.C.~Alperin,  The modular tree of Pythagoras,  Amer. Math. Monthly
\textbf{112} (2005),  807-816.  


\bibitem[BT]{BT}G.~Benkart and P.~Terwilliger,  
The universal central extension  
of the three-point $\mathfrak {sl}_2$ loop algebra,
Proc. Amer. Math. Soc. \textbf{135} (2007), 1659--1668.

\bibitem[H1]{H1} B.~Hartwig,  Tridiagonal Pairs, the Onsager Algebra, and the
Three-Point $\mathfrak{sl}_2$ Loop Algebra, Ph.D. thesis, University of
Wisconsin-Madison, 2006.  
 
 \bibitem[H2]{H2} B.~Hartwig, The tetrahedron algebra and its finite-dimensional irreducible modules,
 Linear Algebra Appl. \textbf{422} (2007),  219--235. 
 
\bibitem[HT]{HT} B.~Hartwig and P.~Terwilliger, 
The Tetrahedron algebra,
the Onsager algebra, and the
$\mathfrak{sl}_2$ loop algebra, 
J. Algebra \textbf{308} (2007), 840--863.

{\tt arXiv:math-ph/0511004}.

\bibitem[Hu]{Hu} J.~Humphreys, \emph{Introduction to Lie Algebras
and Representation Theory}, Springer Verlag, New York, 1972.
 
\bibitem[ITW]{ITW} T.~Ito, P.~Terwilliger, C.~Weng,
The quantum algebra $U_q(\mathfrak {sl}_2)$ and its equitable presentation,
J. Algebra \textbf{208} (2006),  284-301.

\bibitem[K]{K} V.G.~Kac, \emph{Infinite Dimensional Lie Algebras}, 
Third Ed., Cambridge U.
Press, Cambridge, 1990.  

\bibitem[M]{M} R.V.~Moody, Root systems of hyperbolic type,  Adv. in Math. 
\textbf{33} (1979), no. 2, 144--160. 

\bibitem[MP]{MP} R.V.~Moody and A.~Pianzola, \emph{Lie Algebras with
Triangular Decomposition}, Wiley,  New
York 1995.  
%%%%%%%%%%%%%%%%%%%%%%%%%%%%%%%%%%%%%%%%%%%%%%%%%%%%
%\bibitem[R]{R} D.~Romik,  The dynamics of Pythagorean triples, arX.
%%%%%%%%%%%%%%%%%%%%%%%%%%%%%%%%%%%%%%%%%%%%%%%%%%%%%%

%\bibitem[S]{S} R.~Steinberg, Automorphisms of classical Lie algebras, 
% Pacific J. Math \textbf{11} (1961),1119-1129.


\bibitem[TW]{TW} I.~Tuba and H.~Wenzl,  Representations of the braid
group B$_3$ and of SL$(2,\mathbb Z)$,  Pacific J. Math \textbf{197} (2001),
491-510.  


\end{thebibliography}
\end{document}